\documentclass[12pt]{article}
\usepackage{amsmath}
\usepackage{epsf}
\usepackage{amssymb}
\usepackage{graphicx}
\usepackage{ifpdf}
\usepackage{caption}

\usepackage{array}
\usepackage{booktabs} 
\usepackage{multirow}

\newtheorem{theorem}{Theorem}[section]

\newtheorem{definition}{Definition}[section]
\newtheorem{lemma}{Lemma}[section]

\begin{document}
\pagestyle{plain}
\renewcommand{\thefootnote}{\fnsymbol{footnote}}

\begin{center}
{\bf \Large The generalized 4-connectivity of burnt pancake graphs}\footnote{This
work was supported by National Science Foundation of China(No.12271157 \& No.12371346), Natural Science Foundation of Hunan Province (No.2022JJ30028 \& No.2023JJ30072).}
\vskip 5mm

{{\bf Jing Wang$^1$, Jiang Wu$^1$, Zhangdong Ouyang$^2$, Yuanqiu Huang$^3$}\\[2mm]
$^1$ School of Mathematics, Changsha University, Changsha, China\\
$^2$ School of Mathematics, Hunan First Normal University, Changsha, China \\
$^3$ School of Mathematics, Hunan Normal University, Changsha, China}\\[6mm]
\end{center}
\date{}

\noindent{\bf Abstract}\; The generalized $k$-connectivity of a graph $G$, denoted by $\kappa_k(G)$, is the minimum number of internally edge disjoint $S$-trees for any $S\subseteq V(G)$ and $|S|=k$. The generalized $k$-connectivity is a natural extension of the classical connectivity and plays a key role in applications related to the modern interconnection networks. An $n$-dimensional burnt pancake graph $BP_n$ is a Cayley graph which posses many desirable properties. In this paper, we try to evaluate the reliability of $BP_n$ by investigating its generalized 4-connectivity. By introducing the definition of inclusive tree and by studying structural properties of $BP_n$, we show that $\kappa_4(BP_n)=n-1$ for $n\ge 2$, that is, for any four vertices in $BP_n$, there exist ($n-1$) internally edge disjoint trees connecting them in $BP_n$.

{\bf MR(2000) Subject Classification} 05C40, 05C05

\section{Introduction}\label{secintro}

With rapid development and advances of very large scale integration technology and wafer-scale integration technology, multiprocessor systems have been widely designed and used in our daily life. It is well known that the underlying topology of the multiprocessor systems can be modelled by a connected graph $G=(V(G),E(G))$, where $V(G)$ is the set of processors and $E(G)$ is the set of communication links of multiprocessor systems.

Fault tolerance has become increasingly significant nowadays since multiprocessor systems failure is inevitable. The connectivity  is a key parameter for  measuring fault tolerance of the network. A subset $S\subseteq V(G)$ of a connected graph $G$ is called a {\it vertex-cut} if $G\setminus S$ is disconnected or trivial. The {\it connectivity} $\kappa(G)$ of $G$ is defined as the minimum cardinality over all vertex-cuts of $G$. Note that the larger $\kappa(G)$ is, the more reliable the network is. A well known theorem of Whitney \cite{Whitney1932} provides an equivalent definition of connectivity. For each 2-subset $S=\{x,y\}\subseteq V(G)$, let $\kappa(S)$ denote the maximum number of internally disjoint ($x,y$)-paths in $G$. Then
\begin{equation*}
\kappa(G)=\min\{\kappa(S) ~|~S\subseteq V(G)\; {\rm and} \; |S|=2\}.
\end{equation*}

The generalized $k$-connectivity, which was introduced by Chartrand et al. \cite{Chartrand1984}, is a strengthening of connectivity and can be served as an essential parameter for measuring reliability and fault tolerance of the network. Let $G=(V(G),E(G))$ be a simple graph, $S$ be a subset of $V(G)$. A tree $T$ in $G$ is called an $S$-{\it tree}, if $S\subseteq V(T)$. The trees $T_1, T_2, \cdots, T_r$ are called {\it internally edge disjoint $S$-trees} if $V(T_i)\cap V(T_j)=S$ and $E(T_i)\cap E(T_j)=\emptyset$ for any integers $1\le i\ne j\le r$. Let $\kappa_G(S)$ denote the maximum number of internally edge disjoint $S$-trees. For an integer $k$ with $2\le k\le |V(G)|$, the {\it generalized $k$-connectivity} of $G$, denoted by $\kappa_k(G)$, is defined as
\begin{equation*}
\kappa_k(G)=\min\{\kappa_G(S) ~|~S\subseteq V(G)\; {\rm and} \; |S|=k\}.
\end{equation*}

The generalized 2-connectivity is exactly the classical connectivity. Over the past few years, research on the generalized connectivity has received meaningful progress. Li et al. \cite{SLi2012n} derived that it is NP-complete for a general graph $G$ to decide whether there are $l$ internally disjoint trees connecting $S$, where $l$ is a fixed integer and $S\subseteq V(G)$. Authors in \cite{HZLi2014,SLi2010} investigated the upper and lower bounds of the generalized connectivity of a general graph $G$.

Many authors tried to study exact values of the generalized connectivity of graphs. The generalized $k$-connectivity of the complete graph, $\kappa_k(K_n)$, was determined in \cite{Chartrand2010} for every pair $k,n$ of integers with $2\le k\le n$. The generalized $k$-connectivity of the complete bipartite graph $K_{a,b}$ was obtained in \cite{SLi2012b} for all $2\le k\le a+b$. For $k=3$ or $k=4$, the generalized $k$-connectivity of other important classes of graphs, such as, Cartesian product graphs \cite{HZLi2012,HZLi2017}, hypercubes \cite{HZLi2012,SLin2017}, dual cubes \cite{ZhaoHao20191}, exchanged hypercubes \cite{K4EHst}, balanced hypercubes \cite{Wei2021}, locally twisted cubes \cite{K4LQn,Wang2021}, hierarchical cubic networks \cite{K4HSn}, folded hypercubes \cite{FQn}, divide-and-swap cubes \cite{K4DSCn}, star graphs and bubble-sort graphs \cite{SLi2016}, bubble-sort star graphs \cite{Hao20191}, ($n,k$)-star graphs \cite{Snk2020,ANn12018}, pancake graphs \cite{K4Pn}, several Cayley graphs \cite{SLi2017,ZHao20193,ANn12019} et al. have draw many scholars' attention. As we can see, the results on the generalized $k$-connectivity of networks are almost about $k\le 4$.

Gates (the founder of Microsoft) and Papadimitriou \cite{Gate1979} introduced the Burnt Pancake Problem in 1979. With deep understanding, it is known that the Burnt Pancake Problem relates to the construction of networks of parallel processors. Many scholars worked a lot on structural properties of burnt pancake graph $BP_n$. In \cite{Chin2009},  Chin et al. proved  $BP_n$ to be regular and vertex-transitive. Moreover, the spanning connectivity of $BP_n$ was also determined in \cite{Chin2009}. Lai and Yu \cite{LaiYu2011} proved that $BP_n$ contains $n$ independent Hamiltonian cycles. In \cite{Song2015}, Song et al. found that $BP_n$ is a Cayley graph, in addition, they investigated the super and extra connectivity of $BP_n$. In this paper, we try to evaluate the reliability of $BP_n$ by studying its generalized 4-connectivity and obtain the following result.

\begin{theorem}\label{thmBPn}
For $n\ge 2$, $\kappa_4(BP_n)=n-1$.
\end{theorem}

This paper is organized as follows. Section \ref{secpre} introduces some necessary preliminaries. Section \ref{secBPn} presents the definition and structural properties of burnt pancake graph $BP_n$. Let $S$ be any 4-subset of $V(BP_n)$. In order to prove Theorem \ref{thmBPn}, we shall make a lot of efforts in Sections \ref{secS3}, \ref{secS2} and \ref{secS1} to prove that there are ($n-1$)-internally edge disjoint $S$-trees in $BP_n$, depending on the maximum cardinality of $S\cap V(G^i)$ for $i\in [[n]]$. Then Theorem \ref{thmBPn} could be proved in Section \ref{secthm}. Finally, the paper is concluded in Section \ref{seccon}.

\section{Preliminary}\label{secpre}

First of all, we introduce some necessary preliminaries. Let $G=(V(G),E(G))$ be a simple and connected graph with $V(G)$ be its vertex set and $E(G)$ be its edge set.
For a vertex $x\in V(G)$, the {\it degree} of $x$ in $G$, denoted by ${\rm deg}_G(x)$, is the number of edges of $G$ incident with $x$. Denote $\delta(G)$ the {\it minimum degree} of vertices of $G$. A graph is $r$-{\it regular} if ${\rm deg}_G(x)=r$ for every vertex $x\in V(G)$. For a vertex $x\in V(G)$, we use $N_G(x)$ to denote the neighbour set of $x$ and $N_G[x]$ to denote $N_G(x)\cup \{x\}$. Let $V'\subseteq V(G)$, denote $G\backslash V'$ be the graph obtained from $G$ by deleting all the vertices in $V'$ together with their incident edges. Denote by $G[V']$ the subgraph of $G$ induced on $V'$.

Let $P$ be a path in $G$ with $x$ and $y$ be its two terminal vertices, then $P$ is called an $(x,y)$-{\it path}. Two ($x,y$)-paths $P_1$ and $P_2$ are {\it internally disjoint} if they have no common internal vertices, that is, $V(P_1)\cap V(P_2)=\{x, y\}$.

The following Lemma \ref{lemxypath}, Lemma \ref{lemKfan} and Lemma \ref{lemXYpaths} are results on the connectivity of a graph that are well-known in the literature.

\begin{lemma}\label{lemxypath} (\cite{Bondy})
Let $G$ be a $k$-connected graph, and let $x$ and $y$ be a pair of distinct vertices of $G$. Then there exist $k$ internally disjoint ($x,y$)-paths in $G$.
\end{lemma}

\begin{lemma}\label{lemKfan} (\cite{Bondy})
Let $G$ be a $k$-connected graph, let $x$ be a vertex of $G$ and let $Y\subseteq V(G)\backslash\{x\}$ be a set of at least $k$ vertices of $G$. Then there exists a $k$-fan in $G$ from $x$ to $Y$, that is, there exists a family of $k$ internally disjoint ($x,Y$)-paths whose terminal vertices are distinct in $Y$.
\end{lemma}

\begin{lemma}\label{lemXYpaths} (\cite{Bondy})
Let $G$ be a $k$-connected graph, and let $X$ and $Y$ be subsets of $V(G)$ of cardinality at least $k$. Then there exists a family of $k$ pairwise disjoint ($X,Y$)-paths in $G$.
\end{lemma}

For $3\le k\le |V(G)|$, Li gave an upper bound of $\kappa_k(G)$ for a general graph $G$ in her Ph.D. thesis \cite{LiSSPhD2012}.

\begin{lemma}\label{lemupperKk}(\cite{LiSSPhD2012})
Let $G$ be a connected graph with minimum degree $\delta(G)$. If there are two adjacent vertices of degree $\delta(G)$, then $\kappa_k(G)\le \delta(G)-1$ for $3\le k\le |V(G)|$.
\end{lemma}

The following result is about the relationship between $\kappa_k(G)$ and $\kappa_{k-1}(G)$ of a regular graph $G$.

\begin{lemma}\label{lemKkk-1}(\cite{SLin2017})
Let $G$ be an $r$-regular graph. If $\kappa_k(G)=r-1$, then $\kappa_{k-1}(G)=r-1$, where $k\ge 4$.
\end{lemma}

\section{The burnt pancake graph $BP_n$}\label{secBPn}

Let $[n]=\{1,2,\cdots, n\}$. For an integer $i$, it is well known that $|i|$ denotes the absolute value of $i$. Denote $\overline{i}=-i$ in this paper. Let $[[n]]$ be the set $[n]\cup \{\overline{i} ~|~ i\in[n]\}$. A {\it signed permutation} of $[n]$ is an $n$-permutation $x_1x_2\cdots x_n$ of $[[n]]$ such that $|x_1||x_2|\cdots |x_n|$ forms a permutation of $[n]$. For a signed permutation $x=x_1x_2\cdots x_n$ of $[[n]]$ and an integer $i$ ($1\le i\le n$), the $i$th {\it signed prefix reversal} of $x$ is denoted by $x(i)=\overline{x}_i \overline{x_{i-1}}\cdots \overline{x_1}x_{i+1}\cdots x_n$.

\begin{figure}[htbp]
\begin{minipage}[t]{0.28\linewidth}
\centering
\resizebox{0.7\textwidth}{!} {\includegraphics{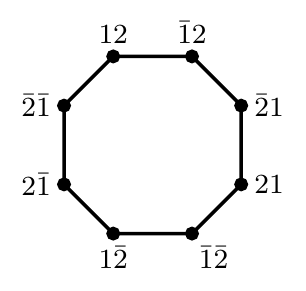}}
\caption{\small The burnt pancake graph $BP_2$} \label{figBP2}
\end{minipage}
\begin{minipage}[t]{0.7\linewidth}
\centering
\resizebox{0.85\textwidth}{!} {\includegraphics{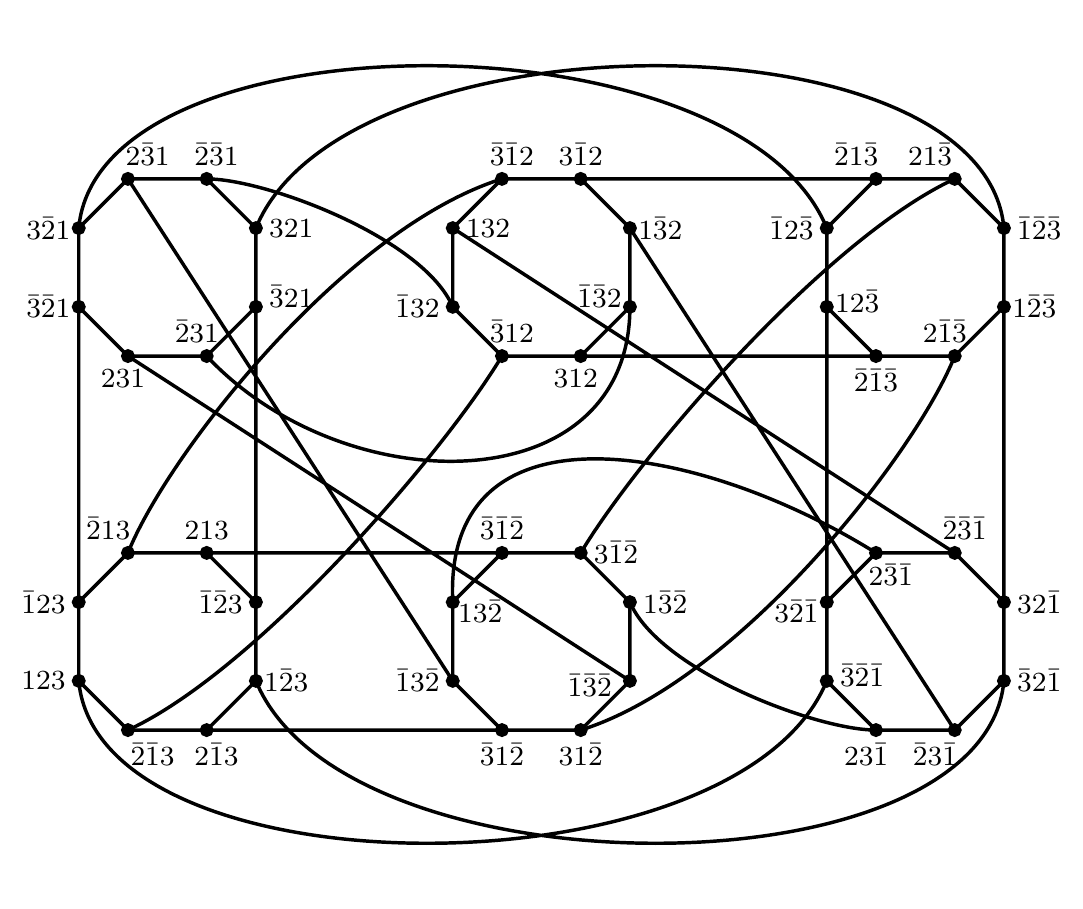}}
\caption{\small The burnt pancake graph $BP_3$} \label{figBP3}
\end{minipage}
\end{figure}

\begin{definition}\label{defBPn}(\cite{Chin2009})
For $n\ge 2$, an $n$-dimensional burnt pancake graph $BP_n$ is a graph with vertex set $V(BP_n)=\{x\, | x$ is a signed permutation of $[[n]]\}$. Two vertices $x=x_1x_2\cdots x_n$ and $y=y_1y_2\cdots y_n$ are adjacent in $BP_n$ if and only if there exists an integer $i$ ($1\le i\le n$) such that $x(i)=y$.
\end{definition}

The burnt pancake graphs  $BP_2$ and $BP_3$ are depicted in Figure \ref{figBP2} and Figure \ref{figBP3}, respectively. Note that, by fixing the symbol $i$ in the rightmost position of each vertex for $i\in[[n]]$, $BP_n$ can be decomposed into $2n$ vertex disjoint subgraphs $BP_n^i$, called {\it clusters}. Obviously, $BP_n^i$ is isomorphic to $BP_{n-1}$ for each $i\in[[n]]$. We write the construction of $BP_n$ symbolically as
$$BP_n=BP_n^1\oplus BP_n^{\overline{1}}\oplus BP_n^2 \oplus BP_n^{\overline{2}}\oplus \cdots \oplus  BP_n^n \oplus BP_n^{\overline{n}}.$$

For simplicity, we shall use $G^i$  to replace $BP_n^i$ in the following texts for $i\in[[n]]$. That means
$$BP_n=G^1\oplus G^{\overline{1}}\oplus G^2 \oplus G^{\overline{2}}\oplus \cdots \oplus  G^n \oplus G^{\overline{n}}.$$

Let $x=x_1x_2\cdots x_n$ be a vertex in $V(G^{x_n})$, it is seen that $x(n)$, the $n$th-signed prefix reversal of $x$, belong to a different cluster $G^{\overline{x_1}}$. We say that  $xx(n)$ is a {\it cross edge} between $G^{x_n}$ and $G^{\overline{x_1}}$. Moreover, we say that $x(n)$ is the {\it out-neighbour} of $x$ and denote it by $\widehat{x}$. The following Lemma \ref{lemBPn1} presents some properties of $BP_n$.

\begin{lemma}\label{lemBPn1}(\cite{Compeau2011,Gate1979,Song2015,BPN3})
For $n\ge 2$, the burnt pancake graph $BP_n$ has the following properties:\\
(1) It is an $n$-regular Cayley graph with $2^nn!$ vertices and $n\times n! \times 2^{n-1}$ edges;\\
(2) $\kappa(BP_n)=n$ and $\kappa_3(BP_n)=n-1$;\\
(3) The girth of $BP_n$ is 8;\\
(4) For $i\in [[n]]$, no two vertices in $V(G^i)$  have a common out-neighbour. Furthermore, let $x=x_1x_2\cdots x_{n}$ be a vertex in $V(G^{x_n})$, the out-neighbours of vertices in $N_{G^{x_n}}[x]$ belong to $n$ different clusters of $BP_n$;\\
(5) For $\{i,j\}\subseteq [[n]]$, we denote by $E(G^i,G^j)$ the set of cross edges between $G^i$ and $G^j$. Then
\[
|E(G^i,G^j)|= \left\{
\begin{array}{ll}
 0, \quad &\;{\rm if}\; i= \overline{j}; \\ \\
(n-2)!\times 2^{n-2}, \quad & \;{\rm if}\;  i\ne \overline{j}.
\end{array}
\right.
\]

\end{lemma}

The following result is not difficult to obtain. Since it appears several times in the paper, we state it formally.

\begin{lemma}\label{lemBPn2}
Let $n\ge 3$ and $BP_n=G^1\oplus G^{\overline{1}}\oplus G^2 \oplus G^{\overline{2}}\oplus \cdots \oplus  G^n \oplus G^{\overline{n}}$. If $x$ is a vertex in $V(G^i)$ with $\widehat{x}\in V(G^j)$, where $i\in [[n]]$ and $j\in [[n]]\backslash\{i,\overline{i}\}$. Then $\widehat{x}(1)\in V(G^{\overline{j}})$.
\end{lemma}

\begin{lemma}\label{lemBPn3}
Let $n\ge 3$ and $BP_n=G^1\oplus G^{\overline{1}}\oplus G^2 \oplus G^{\overline{2}}\oplus \cdots \oplus  G^n \oplus G^{\overline{n}}$. If $x$ is a vertex in $V(G^j)$, where $j\in [[n]]$. Then $G^j\backslash\{x,x(i)\}$ is connected,  $1\le i \le n-1$.
\end{lemma}

\noindent{\bf Proof}\; For $n\ge 4$, $G^j\backslash\{x,x(i)\}$ is connected since $\kappa(G^j)=\kappa(BP_{n-1})=n-1>2$ according to Lemma \ref{lemBPn1}(2). For $n=3$, note that $G^j$ is isomorphic to a cycle of length 8, moreover, $x$ and $x(i)$ are adjacent to each other in $G^j$. Hence, $G^j\backslash\{x,x(i)\}$ is connected for $1\le i \le n-1$. \hfill$\Box$

\vskip 2mm

Let $n\ge 2$ and $x=x_1x_2\cdots x_n$ be a vertex in $V(G^{x_n})$. For $i\in [n]\backslash \{x_n,\overline{x}_n\}$, there is a vertex in $N_{G^{x_n}}(x)$ whose out-neighbour either belong to $V(G^i)$ or belong to $V(G^{\overline{i}})$. For $i\in [n]\backslash \{x_n, \overline{x}_n\}$, denote by $\Gamma_i(x)$ the vertex in $N_{G^{x_n}}(x)$ that its out-neighbour $\widehat{\Gamma}_i(x)\in V(G^i)\cup V(G^{\overline{i}}).$
Throughout this paper, we use $P(x,\widehat{\Gamma}_i(x))$ to denote the path $\{x\Gamma_i(x), \Gamma_i(x)\widehat{\Gamma}_i(x)\}$. That means
$$P(x,\widehat{\Gamma}_i(x))=\{x\Gamma_i(x), \Gamma_i(x)\widehat{\Gamma}_i(x)\}.$$

Let $BP_n=G^1\oplus G^{\overline{1}}\oplus G^2 \oplus G^{\overline{2}}\oplus \cdots \oplus  G^n \oplus G^{\overline{n}}$ and $S=\{x,y,z,w\}$ be any 4-subset of $V(BP_n)$ for $n\ge 3$. Our following discussions in Section \ref{secS3}, Section \ref{secS2} and Section \ref{secS1} are based on the maximum value of $|S\cap V(G^i)|$ for $i\in [[n]]$.

\section{ $\max\{|S\cap V(G^i)|\}=3$ for $i\in [[n]]$}\label{secS3}

For convenience, hereafter we will use IDSTs and IDPs to represent internally edge disjoint $S$-trees and internally disjoint paths, respectively.

\begin{lemma}\label{lemK4S3}
For $n\ge 3$, let $BP_n=G^1\oplus G^{\overline{1}}\oplus G^2 \oplus G^{\overline{2}}\oplus \cdots \oplus  G^n \oplus G^{\overline{n}}$ and  $S=\{x,y,z,w\}$ be any 4-subset of $V(BP_n)$. If there is an integer $i\in [[n]]$ such that $|S\cap V(G^i)|=3$. Then there exist ($n-1$) IDSTs in $BP_n$.
\end{lemma}

\noindent{\bf Proof}\;  Without loss of generality, we may assume that $\{x,y,z\}\subseteq V(G^n)$ and $w\in V(G^j)$, where $j\in [[n]]\backslash \{n\}$. Furthermore, we may assume that $x=12\cdots n$. Note that  $\Gamma_i(x)=x(i)$ and $\widehat{\Gamma}_i(x)\in V(G^i)$ for $i\in [n-1]$.

Let
\begin{equation}\label{eqI1}
I_1=\{i\in [n-1] ~|~ \{\widehat{\Gamma}_i(x), \widehat{\Gamma}_i(z)\}\subseteq V(G^i) ~{\rm and} ~\widehat{\Gamma}_i(y)\in V(G^{\overline{i}})\},
\end{equation}
\begin{equation}\label{eqI2}
I_2=\{i\in [n-1] ~|~ \widehat{\Gamma}_i(x)\in V(G^i) ~{\rm and} ~ \{\widehat{\Gamma}_i(y), \widehat{\Gamma}_i(z)\}\subseteq V(G^{\overline{i}})\},
\end{equation}
\begin{equation}\label{eqI3}
I_3=\{i\in [n-1] ~|~ \{\widehat{\Gamma}_i(x), \widehat{\Gamma}_i(y)\}\subseteq V(G^i) ~{\rm and} ~\widehat{\Gamma}_i(z)\in V(G^{\overline{i}})\},
\end{equation}
\begin{equation}\label{eqI4}
I_4=\{i\in [n-1] ~|~ \{\widehat{\Gamma}_i(x), \widehat{\Gamma}_i(y), \widehat{\Gamma}_i(z)\}\subseteq V(G^i)\}.
\end{equation}

By Lemma \ref{lemBPn1}(3), the subgraph $G^n[\{x,y,z\}]$ contains at most two edges.

\vskip 2mm
{\bf Case 1}. The subgraph $G^n[\{x,y,z\}]$ is empty.

For simplicity, we may assume that $I_1=\{1,\cdots, l_1\}$, $I_2=\{l_1+1,\cdots, l_2\}$, $I_3=\{l_2+1,\cdots,l_3\}$ and $I_4=\{l_3+1,\cdots, n-1\}$. It is possible that $I_i=\emptyset$ for $i\in [4]$. The possibility that $I_i=\emptyset$ for $i\in [3]$ will not affect the following discussions.

\vskip 2mm
{\bf Subcase 1.1}. $j=\overline{n}$.

That means $w\in V(G^{\overline{n}})$. By Lemma \ref{lemBPn1}(5), there are ($n-1$) vertices $w_1, \cdots, w_{n-1}$ in $V(G^{\overline{n}})$ that $\widehat{w}_i\in V(G^i)$ for $i\in [n-1]$. According to Lemma \ref{lemKfan} and the fact that $\kappa(G^{\overline{n}})=\kappa(BP_{n-1})=n-1$, there is an ($n-1$)-fan $P_1, \cdots, P_{n-1}$ in $G^{\overline{n}}$ from $w$ to $\{w_1, \cdots, w_{n-1}\}$ that $w_i\in V(P_i)$ for $i\in [n-1]$.

\vskip 2mm
{\bf Subcase 1.1.1}. $I_4\ne \emptyset$.

That means $(n-1)\in I_4$. For $i\in I_4$, there is a $\{\widehat{\Gamma}_i(x), \widehat{\Gamma}_i(y), \widehat{\Gamma}_i(z), \widehat{w}_i\}$-tree $\widehat{T}_i$ in $G^i$ since $G^i$ is connected. For $i\in I_4$, let
\begin{eqnarray}\label{eqTi5}
T_i=\widehat{T}_i \cup P(x,\widehat{\Gamma}_i(x))\cup P(y,\widehat{\Gamma}_i(y))\cup P(z,\widehat{\Gamma}_i(z))\cup P_i\cup \{w_i\widehat{w}_i\}.
\end{eqnarray}

By Lemma \ref{lemBPn1}(5) and Lemma \ref{lemBPn2}, for $1\le i\le l_3$, there are vertices $a_i$ and $b_i=a_i(1)$ in $V(G^{\overline{n-1}})$ such that $\widehat{a}_i\in V(G^i)$ and  $\widehat{b}_i\in V(G^{\overline{i}})$.

\begin{figure}[htbp]
\centering
\resizebox{0.8\textwidth}{!} {\includegraphics{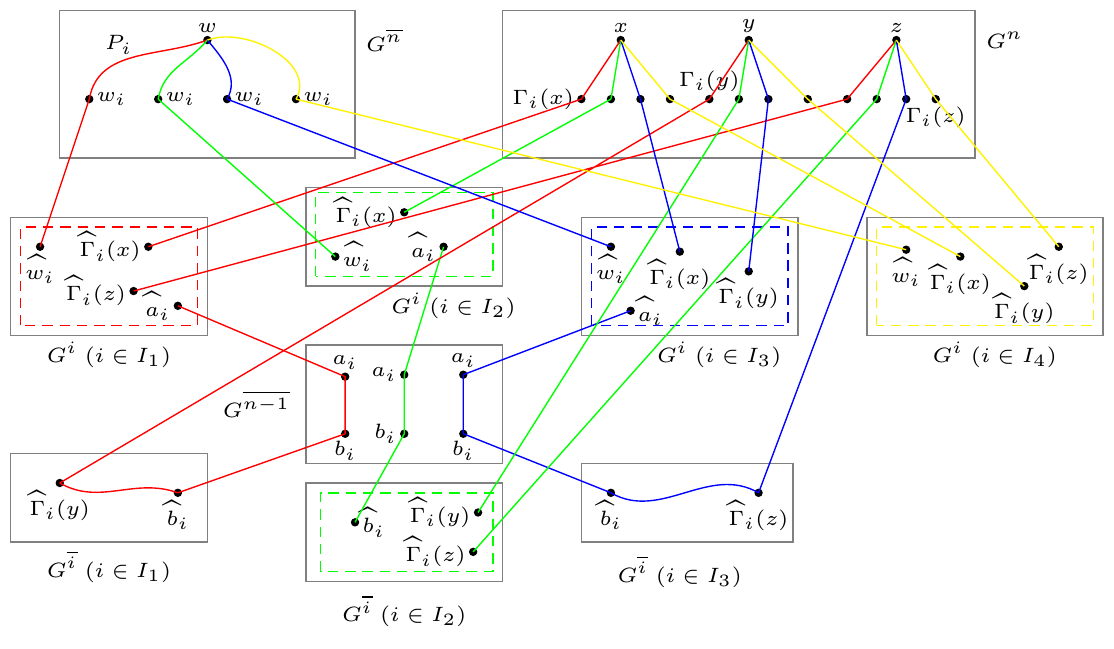}}
\caption{\small Illustration for Subcase 1.1.1 } \label{figS3C11}
\end{figure}

For $i\in I_2$, there is a $\{\widehat{\Gamma}_i(x), \widehat{w}_i, \widehat{a}_i\}$-tree $T'_i$ in $G^i$ since $G^i$ is connected. Moreover, there is a $\{\widehat{b}_i, \widehat{\Gamma}_i(y),\widehat{\Gamma}_i(z)\}$-tree $T''_i$ in $G^{\overline{i}}$. Let
$$IT_i=T'_i\cup T''_i\cup \{a_i\widehat{a}_i, a_ib_i, b_i\widehat{b}_i\}, \;\; i\in I_2.$$
For $i\in I_2$, we say that $IT_i$ is a $\{\widehat{\Gamma}_i(x), \widehat{\Gamma}_i(y), \widehat{\Gamma}_i(z), \widehat{w}_i\}${\it-inclusive tree connects $G^i$ and $G^{\overline{i}}$ passing through $a_ib_i$}. Similarly, there exists a $\{\widehat{\Gamma}_i(x), \widehat{\Gamma}_i(y), \widehat{\Gamma}_i(z), \widehat{w}_i\}$-inclusive tree $IT_i$ connects $G^i$ and $G^{\overline{i}}$ passing through $a_ib_i$ for all $i\in I_1\cup I_3$.

For $1\le i\le l_3$, let
\begin{eqnarray}\label{eqTi6}
T_i=IT_i\cup P(x,\widehat{\Gamma}_i(x))\cup P(y,\widehat{\Gamma}_i(y))\cup P(z,\widehat{\Gamma}_i(z)) \cup P_i\cup \{w_i\widehat{w}_i\}.
\end{eqnarray}
See Figure \ref{figS3C11}. Then $T_1, \cdots, T_{n-1}$ are ($n-1$) IDSTs in $BP_n$.

\vskip 2mm
{\bf Subcase 1.1.2}. $I_4=\emptyset$.

\vskip 2mm
{\bf Subcase 1.1.2.1}. $\min\{|I_1|, |I_2|, |I_3|\}\ge 2$.

By Lemma \ref{lemBPn1}(5), there is a $\{\widehat{\Gamma}_i(x), \widehat{\Gamma}_i(z), \widehat{w}_i, \widehat{\Gamma}_{i+1}(y)\}$-tree $\widehat{T}_i$ in $G^i\cup G^{\overline{i+1}}$ for $1\le i\le l_1-1$. Meanwhile, there is a $\{\widehat{\Gamma}_{l_1}(x), \widehat{w}_{l_1}, \widehat{\Gamma}_{l_1}(z), \widehat{\Gamma}_{1}(y)\}$-tree $\widehat{T}_{l_1}$ in $G^{l_1}\cup G^{\overline{1}}$.

For $1\le i\le {l_1}$, let
$$T_i=\widehat{T}_i\cup  P(x,\widehat{\Gamma}_i(x))\cup P(y,\widehat{\Gamma}_{i+1}(y))\cup P(z,\widehat{\Gamma}_{i}(z)) \cup P_i\cup\{w_i\widehat{w}_i\}.$$
The subscripts are read modulo $l_1$. We may construct other IDSTs $T_{l_1+1}, \cdots, T_{n-1}$ by similar analysis. Therefore, $T_1,$ $\cdots,T_{n-1}$ are ($n-1$) IDSTs in $BP_n$.

\vskip 2mm
{\bf Subcase 1.1.2.2}. $|I_i|\le 1$ for $1\le i\le 3$.

To avoid duplication, we consider the case that $|I_1|=1$, i.e., $I_1=\{1\}$.

Under this situation, it is impossible that $I_2=I_3=\emptyset$ since $n\ge 3$. Without loss of generality, we may assume that $I_2\ne \emptyset$. That means $2\in I_2$. By Lemma \ref{lemBPn1}(5) and Lemma \ref{lemBPn2}, there exists a vertex $u$ in $V(G^{2})\backslash \{\widehat{\Gamma}_2(x),\widehat{w}_2\}$ such that $\widehat{u}\in V(G^1)$ and $\widehat{u}(1)\in V(G^{\overline{1}})$.

The following two subcases are distinguished to construct $T_i$s for $i\in I_1\cup I_2$.

\vskip 2mm
{\bf Subcase 1.1.2.2.1}. $|I_2|=1$.

There exists a vertex $v\in V(G^{\overline{1}})\backslash \{\widehat{\Gamma}_1(y), \widehat{u}(1)\}$ such that $\widehat{v}\in V(G^2)$ and $\widehat{v}(1)\in V(G^{\overline{2}})$. This could be done since $(n-2)!\times 2^{n-2}\ge 2$ when $n\ge 3$. By Lemma \ref{lemBPn3}, both $G^2\backslash \{u, u(1)\}$ and $G^{\overline{1}}\backslash \{v,v(1)\}$ are connected. By similar analysis in Subcase 1.1.1, there is a $\{\widehat{\Gamma}_1(x), \widehat{\Gamma}_1(y), \widehat{\Gamma}_{1}(z), \widehat{w}_1\}$-inclusive tree $IT_1$ connects $G^1$ and $G^{\overline{1}}\backslash \{v,v(1)\}$ passing through $uu(1)$ and a $\{\widehat{\Gamma}_2(x), \widehat{\Gamma}_2(y), \widehat{\Gamma}_{2}(z), \widehat{w}_2\}$-inclusive tree $IT_2$ connects $G^2\backslash \{u,u(1)\}$ and $G^{\overline{2}}$ passing through $vv(1)$, respectively. For $i\in I_1\cup I_2$, let $T_i$ be the same as in Eq.(\ref{eqTi6}).

\vskip 2mm
{\bf Subcase 1.1.2.2.2}. $|I_2|\ge 2$.

By similar arguments in Subcase 1.1.1, there is a $\{\widehat{\Gamma}_1(x), \widehat{\Gamma}_1(y), \widehat{\Gamma}_{1}(z), \widehat{w}_1\}$-inclusive tree $IT_1$ connects $G^1$ and $G^{\overline{1}}$ passing through $uu(1)$. For $i=1$, let $T_i$ be the same as in Eq.(\ref{eqTi6}).

According to Lemma \ref{lemBPn3} and Lemma \ref{lemBPn1}, $\big(G^2\backslash \{u,u(1)\}\big)\cup G^{\overline{3}}$ is connected. Thus, there is a $\{\widehat{\Gamma}_{2}(x), \widehat{w}_2,  \widehat{\Gamma}_{3}(y), \widehat{\Gamma}_{3}(z)\}$-tree $\widehat{T}_2$ in $\big(G^2\backslash \{u,u(1)\}\big)\cup G^{\overline{3}}$. For $i\in I_2$, let $T_i$ be the same as in Subcase 1.1.2.1.

\vskip 2mm

For $i\in I_3$, we may let $T_i$ be the same as in Subcase 1.1.2.1 if $|I_3|\ge 2$ and construct $T_i$ by similar methods in Subcase 1.1.2.2.1 if $|I_3|=1$.

Then $T_1, \cdots,  T_{n-1}$ are ($n-1$) IDSTs in $BP_n$.

\vskip 2mm
{\bf Subcase 1.2}. $\overline{j}\in I_4$.

Without loss of generality, we may assume that $j=\overline{n-1}$. That is to say, $w\in V(G^{\overline{n-1}})$. Remind that $\{\widehat{\Gamma}_i(x), \widehat{\Gamma}_i(y), \widehat{\Gamma}_i(z)\}\subseteq V(G^i)$ for $i\in I_4$. Moreover, either $\widehat{\Gamma}_i(w)\in V(G^i)$ or $\widehat{\Gamma}_i(w)\in V(G^{\overline{i}})$ for $1\le i\le n-2$.

By Lemma \ref{lemBPn2}, there are $2(n-1)$ vertices $a_i$ and $b_i=a_i(1)$ in $V(G^{\overline{n}})$ satisfying $\widehat{a}_i\in V(G^i)$ and $\widehat{b}_i\in V(G^{\overline{i}})$ for $i\in [n-1]$. According to similar analysis in Subcase 1.1.1, there exists a  $\{\widehat{\Gamma}_i(x), \widehat{\Gamma}_i(y), \widehat{\Gamma}_i(z), \widehat{\Gamma}_i(w)\}$-inclusive tree $IT_i$ connects $G^i$ and $G^{\overline{i}}$ passing through $a_ib_i$ for $i\in [n-2]$.

Let $W=\{\Gamma_1(w), \cdots, \Gamma_{n-2}(w)\}.$
Note that $|W|=n-2< \kappa(G^{\overline{n-1}})$, $G^{\overline{n-1}}\backslash W$ is connected. There is a ($w, \widehat{b}_{n-1}$)-path $R_{n-1}$ in $G^{\overline{n-1}}\backslash W$. Furthermore, there is a $\{\widehat{\Gamma}_{n-1}(x), \widehat{\Gamma}_{n-1}(y), \widehat{\Gamma}_{n-1}(z),$ $\widehat{a}_{n-1}\}$-tree $\widehat{T}_{n-1}$ in $G^{n-1}$. Let
\begin{eqnarray*}
T_{n-1}&=&\widehat{T}_{n-1}\cup R_{n-1}\cup \{a_{n-1}\widehat{a}_{n-1}, a_{n-1}b_{n-1}, b_{n-1}\widehat{b}_{n-1}\} \\
&&\cup P(x,\widehat{\Gamma}_{n-1}(x))\cup P(y,\widehat{\Gamma}_{n-1}(y))\cup P(z,\widehat{\Gamma}_{n-1}(z)).
\end{eqnarray*}

For $1\le i\le n-2$, let
\begin{eqnarray}\label{eqTi7}
T_i=IT_i\cup P(x,\widehat{\Gamma}_i(x))\cup P(y,\widehat{\Gamma}_i(y))\cup P(z,\widehat{\Gamma}_{i}(z))\cup P(w,\widehat{\Gamma}_{i}(w)).
\end{eqnarray}
Then $T_1, \cdots,  T_{n-1}$ are ($n-1$) IDSTs in $BP_n$.

\vskip 2mm
{\bf Subcase 1.3}. Either $j\in I_1\cup I_2\cup I_3\cup I_4$ or $\overline{j}\in I_1\cup I_2\cup I_3$.

Without loss of generality, we may assume that $j=n-1$, i.e., $w\in V(G^{n-1})$. For $i\in [n-2]$, define $a_i$ and $b_i$ be the same as in Subcase 1.2. There is a $\{\widehat{\Gamma}_i(x), \widehat{\Gamma}_i(y), \widehat{\Gamma}_i(z), \widehat{\Gamma}_i(w)\}$-inclusive tree $IT_i$ connects $G^i$ and $G^{\overline{i}}$ passing through $a_ib_i$ for $i\in [n-2]$. Furthermore, there is a $\{\widehat{\Gamma}_{n-1}(x), \widehat{\Gamma}_{n-1}(y), \widehat{\Gamma}_{n-1}(z), w\}$-tree $\widehat{T}_{n-1}$ in $G^{n-1}\backslash \{\Gamma_1(w), \cdots, \Gamma_{n-2}(w)\}$.
Let
$$T_{n-1}=\widehat{T}_{n-1}\cup P(x,\widehat{\Gamma}_{n-1}(x))\cup P(y,\widehat{\Gamma}_{n-1}(y))\cup P(z,\widehat{\Gamma}_{n-1}(z))$$
and let $T_i$ ($i\in [n-2]$) be the same as in Eq.(\ref{eqTi7}). Then $T_1, \cdots,  T_{n-1}$ are ($n-1$) IDSTs in $BP_n$.

\vskip 2mm
{\bf Case 2}. The subgraph $G^n[\{x,y,z\}]$ contains exactly one edge.

Without loss of generality, we may assume that $y=x(l)$, where $1\le l\le n-1$. It yields that $y=\overline{l}\cdots \bar{1}(l+1)\cdots n$. Moreover, $\widehat{\Gamma}_i(y)\in V(G^{\overline{i}})$ for $1\le i\le l$ and $\widehat{\Gamma}_i(y)\in V(G^{i})$ for $l+1\le i\le n-1$. Note that $\widehat{\Gamma}_i(z)\in V(G^{i})\cup V(G^{\overline{i}})$ for $i\in [n-1]$.

\vskip 2mm

{\bf Subcase 2.1}. $j=\overline{n}$, i.e., $w\in V(G^{\overline{n}})$.

\vskip 2mm

{\bf Subcase 2.1.1}. $l\ne 1$.

Consider the locations of $\widehat{\Gamma}_1(z)$ and $\widehat{\Gamma}_l(z)$. By former analysis, there are four possibilities: (1) $\widehat{\Gamma}_1(z)\in V(G^1)$ and $\widehat{\Gamma}_l(z)\in V(G^l)$; (2) $\widehat{\Gamma}_1(z)\in V(G^1)$ and $\widehat{\Gamma}_l(z)\in V(G^{\overline{l}})$; (3) $\widehat{\Gamma}_1(z)\in V(G^{\overline{1}})$ and $\widehat{\Gamma}_l(z)\in V(G^l)$;  (4) $\widehat{\Gamma}_1(z)\in V(G^{\overline{1}})$ and $\widehat{\Gamma}_l(z)\in V(G^{\overline{l}})$. By symmetry, we only need to consider the following two cases.

\vskip 2mm

{\bf Subcase 2.1.1.1}.  $\widehat{\Gamma}_1(z)\in V(G^1)$ and $\widehat{\Gamma}_l(z)\in V(G^l)$.

By Lemma \ref{lemBPn1}(5), there are ($n-1$) vertices $w_1, \cdots, w_{n-1}$ in $V(G^{\overline{n}})$ such that $\widehat{w}_i\in V(G^i)$ for $i\in [n-1]$. According to Lemma \ref{lemKfan}, there is a family of ($n-1$) IDPs $P_1, \cdots, P_{n-1}$ in $G^{\overline{n}}$ from $w$ to $\{w_1, \cdots, w_{n-1}\}$  such that $w_i\in V(P_i)$ for $i\in [n-1]$.

Since $G^1\cup G^{\overline{l}}$ is connected, there is a $\{\widehat{\Gamma}_1(x), \widehat{\Gamma}_1(z), \widehat{w}_1, \widehat{\Gamma}_l(y)\}$-tree $\widehat{T}_1$ in $G^1\cup G^{\overline{l}}$. Furthermore, there is a $\{\widehat{y}, \widehat{\Gamma}_l(z), \widehat{w}_l\}$-tree $\widehat{T}_l$ in $G^l$. Let
\begin{eqnarray*}
T_1=\widehat{T}_1\cup  P(x,\widehat{\Gamma}_1(x))\cup P(z,\widehat{\Gamma}_1(z))\cup P(y,\widehat{\Gamma}_l(y))\cup P_1\cup\{w_1\widehat{w}_1\}
\end{eqnarray*}
and
\begin{eqnarray*}
T_l=\widehat{T}_l\cup  P(z,\widehat{\Gamma}_l(z))\cup P_l\cup\{xy, y\widehat{y},w_l\widehat{w}_l\}.
\end{eqnarray*}

\begin{figure}[htbp]
\begin{minipage}[t]{0.49\linewidth}
\centering
\resizebox{0.9\textwidth}{!} {\includegraphics{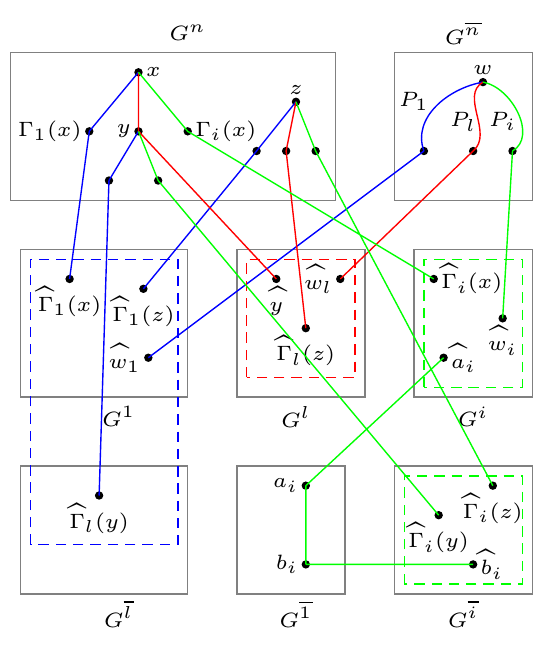}}
\caption{\small Illustration for Subcase 2.1.1.1} \label{figS3C21}
\end{minipage}
\begin{minipage}[t]{0.49\linewidth}
\centering
\resizebox{0.9\textwidth}{!} {\includegraphics{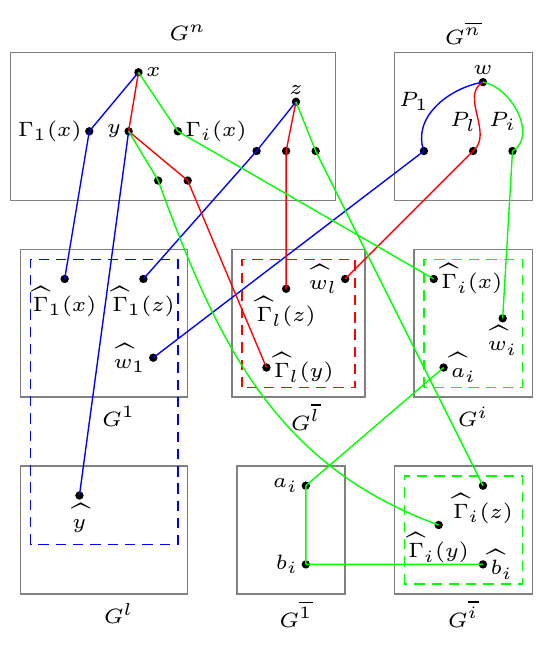}}
\caption{\small Illustration for Subcase 2.1.1.2} \label{figS3C22}
\end{minipage}
\end{figure}

By Lemma \ref{lemBPn2}, for $i\in [n-1]\backslash \{1,l\}$, there are 2($n-3$) vertices $a_i$ and $b_i=a_i(1)$ in $V(G^{\overline{1}})$ such that $\widehat{a}_i\in V(G^i)$ and $\widehat{b}_i\in V(G^{\overline{i}})$. By similar analysis in Subcase 1.1.1, we can find a $\{\widehat{\Gamma}_i(x), \widehat{\Gamma}_i(y), \widehat{\Gamma}_i(z), \widehat{w}_i\}$-inclusive tree $IT_i$ connects $G^i$ and $G^{\overline{i}}$ passing through $a_ib_i$ and construct $T_i$ be the same as in Eq.(\ref{eqTi6}), for $i\in [n-1]\backslash \{1,l\}$. See Figure \ref{figS3C21}. Then $T_1,\cdots, T_{n-1}$ are ($n-1$) IDSTs in $BP_n$.

\vskip 2mm

{\bf Subcase 2.1.1.2}. $\widehat{\Gamma}_1(z)\in V(G^1)$ and $\widehat{\Gamma}_l(z)\in V(G^{\overline{l}})$.

By Lemma \ref{lemBPn1}(5), there are ($n-1$) vertices $w_1, \cdots, w_{n-1}$ in $V(G^{\overline{n}})$ such that $\widehat{w}_l\in V(G^{\overline{l}})$ and  $\widehat{w}_i\in V(G^i)$ for $i\in [n-1]\backslash\{l\}$. According to Lemma \ref{lemKfan}, there is a family of ($n-1$) IDPs $P_1, \cdots, P_{n-1}$ in $G^{\overline{n}}$ from $w$ to $\{w_1, \cdots, w_{n-1}\}$ such that $w_i\in V(P_i)$ for $i\in [n-1]$.

Based on Lemma \ref{lemBPn1}, there is a $\{\widehat{\Gamma}_1(x), \widehat{\Gamma}_1(z), \widehat{w}_1, \widehat{y}\}$-tree $\widehat{T}_1$ in $G^1\cup G^l$ and a $\{\widehat{\Gamma}_l(y), \widehat{\Gamma}_l(z), \widehat{w}_l\}$-tree $\widehat{T}_l$ in $G^{\overline{l}}$. Let
\begin{eqnarray*}
T_1=\widehat{T}_1\cup  P(x,\widehat{\Gamma}_1(x))\cup P(z,\widehat{\Gamma}_1(z))\cup P_1\cup \{y\widehat{y}, w_1\widehat{w}_1\}
\end{eqnarray*}
and
\begin{eqnarray*}
T_l=\widehat{T}_l\cup P(y,\widehat{\Gamma}_l(y))\cup P(z,\widehat{\Gamma}_l(z))\cup  P_l\cup \{xy, w_l\widehat{w}_l\}.
\end{eqnarray*}

For $i\in [n-1]\backslash \{1,l\}$, $IT_i$ and $T_i$ may be defined by similar methods in Subcase 2.1.1.1.  See Figure \ref{figS3C22}. Then $T_1,\cdots, T_{n-1}$ are ($n-1$) IDSTs in $BP_n$.

\vskip 2mm

{\bf Subcase 2.1.2}. $l=1$.

Note that either $\widehat{\Gamma}_1(z)\in V(G^1)$ or $\widehat{\Gamma}_1(z)\in V(G^{\overline{1}})$. By symmetry, we only need to consider the possibility that $\widehat{\Gamma}_1(z)\in V(G^1)$.

There are ($n-1$) vertices $w_1, \cdots, w_{n-1}$ in $V(G^{\overline{n}})$ such that $\widehat{w}_i\in V(G^i)$ for $i\in [n-1]$. According to Lemma \ref{lemKfan}, there is a family of ($n-1$) IDPs $P_1, \cdots, P_{n-1}$ in $G^{\overline{n}}$ from $w$ to $\{w_1, \cdots, w_{n-1}\}$  such that $w_i\in V(P_i)$ for $i\in [n-1]$.

\begin{figure}[htbp]
\begin{minipage}[t]{0.4\linewidth}
\centering
\resizebox{0.9\textwidth}{!} {\includegraphics{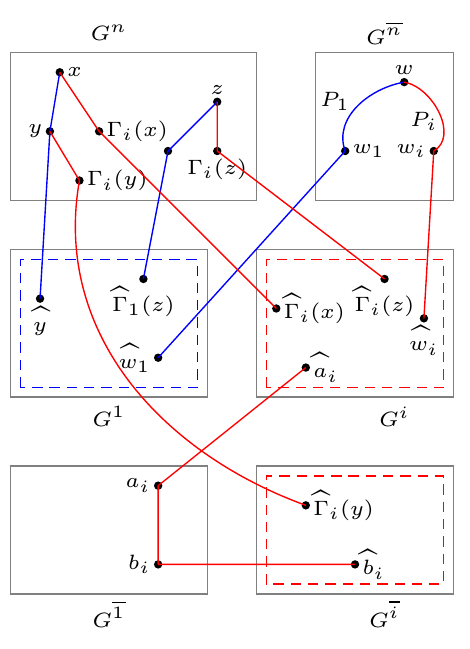}}
\caption{\small Illustration for Subcase 2.1.2} \label{figS3C23}
\end{minipage}
\begin{minipage}[t]{0.55\linewidth}
\centering
\resizebox{0.9\textwidth}{!} {\includegraphics{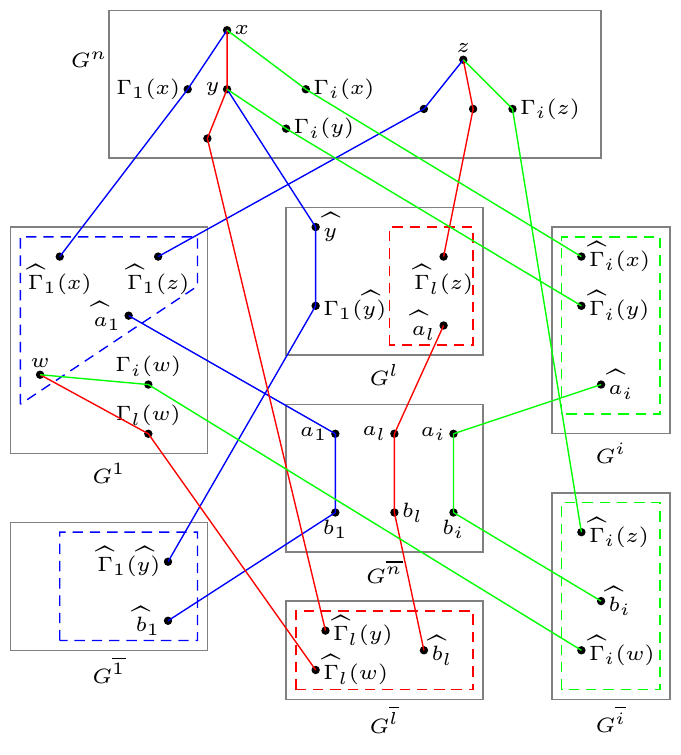}}
\caption{\small Illustration for Subcase 2.2} \label{figS3C24}
\end{minipage}
\end{figure}

Since $G^1$ is connected, there exists a $\{\widehat{y}, \widehat{\Gamma}_1(z), \widehat{w}_1\}$-tree $\widehat{T}_1$ in $G^1$. Let
\begin{eqnarray*}
T_1=\widehat{T}_1\cup  P(z,\widehat{\Gamma}_1(z))\cup P_1\cup\{xy, y\widehat{y}, w_1\widehat{w}_1\}
\end{eqnarray*}
and let $T_i$ be the same as in Subcase 2.1.1.1 for $2\le i \le n-1$. See Figure \ref{figS3C23}. Then $T_1,\cdots, T_{n-1}$ are ($n-1$) IDSTs in $BP_n$.

\vskip 2mm

{\bf Subcase 2.2}. $j\in \{1, \overline{1}, l, \overline{l}\}$.

To avoid duplication, we only consider the case that $j=1$.

Denote by $v=\widehat{\Gamma}_l(y)$ for simplicity. Remind that $\widehat{\Gamma}_1(\widehat{y})\in V(G^1)\cup V(G^{\overline{1}})$ and $\widehat{\Gamma}_1(v)\in V(G^1)\cup V(G^{\overline{1}})$. Since each vertex has an unique out-neighbour, either $\widehat{\Gamma}_1(\widehat{y})\ne \Gamma_l(w)$ or $\widehat{\Gamma}_1(v)\ne \Gamma_l(w)$. Without loss of generality, we may assume that $\widehat{\Gamma}_1(\widehat{y})\ne \Gamma_l(w)$.

There exist 2($n-1$) vertices $a_i$ and $b_i=a_i(1)$ in $V(G^{\overline{n}})$ such that $\widehat{a}_i\in V(G^i)$ and $\widehat{b}_i\in V(G^{\overline{i}})$ for $i\in [n-1]$ according to Lemma \ref{lemBPn2}.

Let $W=\{\Gamma_2(w), \cdots, \Gamma_{n-1}(w)\}$. Note that $|W|=n-2<\kappa(G^1)$. Then $G^1\backslash W$ is connected. By similar arguments in Subcase 1.1.1, there is a $\{\widehat{\Gamma}_1(x), w, \widehat{\Gamma}_1(z), \widehat{\Gamma}_1(\widehat{y})\}$-inclusive tree $IT_1$ connects $\big(G^1\backslash W\big)$ and $G^{\overline{1}}$ passing through $a_1b_1$. Let
\begin{eqnarray*}
T_1=IT_1\cup  P(x,\widehat{\Gamma}_1(x))\cup P(z,\widehat{\Gamma}_1(z))\cup \{y\widehat{y}, ~\widehat{y}\Gamma_1(\widehat{y}), ~\Gamma_1(\widehat{y})\widehat{\Gamma}_1(\widehat{y})\}.
\end{eqnarray*}

It is seen that $G^l\backslash \{\widehat{y}, \Gamma_1(\widehat{y})\}$ is connected by Lemma \ref{lemBPn3}. Using similar analysis in Subcase 1.1.1, there is a $\{\widehat{\Gamma}_l(y), \widehat{\Gamma}_l(z), \widehat{\Gamma}_l(w)\}$-inclusive tree $IT_l$ connects $\big(G^l\backslash \{\widehat{y}, \Gamma_1(\widehat{y})\}\big)$ and $G^{\overline{l}}$ passing through $a_lb_l$. Let
\begin{eqnarray*}
T_l=IT_l\cup P(y,\widehat{\Gamma}_l(y))\cup P(z,\widehat{\Gamma}_l(z))\cup  P(w,\widehat{\Gamma}_l(w))\cup \{xy\}.
\end{eqnarray*}

Again, by similar analysis in Subcase 2.1.1, there is a $\{\widehat{\Gamma}_i(x), \widehat{\Gamma}_i(y), \widehat{\Gamma}_i(z), \widehat{\Gamma}_i(w)\}$-inclusive tree $IT_i$ connects $G^i$ and $G^{\overline{i}}$ passing through $a_ib_i$ for $i\in [n-1]\backslash\{1,l\}$. For $i\in [n-1]\backslash\{1,l\}$, let $T_i$ be the same as in Eq.(\ref{eqTi7}). See Figure \ref{figS3C24}. Then $T_1,\cdots, T_{n-1}$ are ($n-1$) IDSTs in $BP_n$.

\vskip 2mm

{\bf Subcase 2.3}. Either $j\in [n-1]\backslash\{1,l\}$ or $\overline{j}\in [n-1]\backslash\{1,l\}$ .

For simplicity, we may assume that $j=n-1$.

For $i\in [n-1]$, let $a_i$ and $b_i=a_i(1)$ be the same as in Subcase 2.2. In addition, let $W=\{\Gamma_1(w), \cdots, \Gamma_{n-2}(w)\}$.

Since $G^{n-1}\backslash W$ is connected, there is a $\{w, \widehat{\Gamma}_{n-1}(x), \widehat{\Gamma}_{n-1}(y), \widehat{\Gamma}_{n-1}(z)\}$-inclusive tree $IT_{n-1}$ connects $\big(G^{n-1}\backslash W\big)$ and $G^{\overline{n-1}}$ passing through $a_{n-1}b_{n-1}$. Let
\begin{eqnarray*}
T_{n-1}=IT_{n-1}\cup  P(x,\widehat{\Gamma}_{n-1}(x))\cup P(y,\widehat{\Gamma}_{n-1}(y))\cup P(z,\widehat{\Gamma}_{n-1}(z)).
\end{eqnarray*}

\vskip 2mm

{\bf Subcase 2.3.1}. $l\ne 1$.

There is a $\{\widehat{\Gamma}_1(x), \widehat{\Gamma}_1(\widehat{y}), \widehat{\Gamma}_1(z), \widehat{\Gamma}_1(w)\}$-inclusive tree $IT_1$ connects $G^1$ and $G^{\overline{1}}$ passing through $a_1b_1$ since $\{\widehat{\Gamma}_1(\widehat{y}), \widehat{\Gamma}_1(w)\}\subseteq V(G^1)\cup V(G^{\overline{1}})$.  Furthermore, there is a $\{\widehat{\Gamma}_l(y), \widehat{\Gamma}_l(z), \widehat{\Gamma}_l(w)\}$-inclusive tree $IT_l$ connects $(G^l\backslash\{\widehat{y}, \Gamma_1(\widehat{y})\})$ and $G^{\overline{l}}$ passing through $a_lb_l$.

Let
\begin{eqnarray*}
T_{1}=IT_{1}\cup  P(x,\widehat{\Gamma}_{1}(x))\cup P(z,\widehat{\Gamma}_{1}(z))\cup P(w,\widehat{\Gamma}_{1}(w))\cup P(\widehat{y},\widehat{\Gamma}_1(\widehat{y}))\cup \{y\widehat{y}\}.
\end{eqnarray*}

For $2\le i\le n-2$, we may define $T_i$ to be the same as in Subcase 2.2. Then $T_1,\cdots, T_{n-1}$ are ($n-1$) IDSTs in $BP_n$.

\vskip 2mm

{\bf Subcase 2.3.2}. $l=1$.

There is a $\{\widehat{y}, \widehat{\Gamma}_{1}(z), \widehat{\Gamma}_{1}(w)\}$-inclusive tree $IT_1$ connects $G^1$ and $G^{\overline{1}}$ passing through $a_1b_1$. Let
\begin{eqnarray*}
T_{1}=IT_{1}\cup P(z,\widehat{\Gamma}_{1}(z))\cup P(w,\widehat{\Gamma}_{1}(w))\cup \{xy, y\widehat{y}\}.
\end{eqnarray*}

For $2\le i\le n-2$, we may define $T_i$ to be the same as in Subcase 2.2. Then $T_1,\cdots, T_{n-1}$ are ($n-1$) IDSTs in $BP_n$.

\vskip 2mm
{\bf Case 3}. The subgraph $G^n[\{x,y,z\}]$ contains two edges.

Without loss of generality, we may assume that $\{xy, xz\}\subseteq E(G^n)$. Moreover, let $y=x(l)=\overline{l}\cdots \bar{1}(l+1)\cdots n$ and $z=x(k)=\overline{k}\cdots \overline{2}\bar{1}(k+1)\cdots n$, where $1\le l<k\le n-1$. It is seen that $\Gamma_1(y)=x=\Gamma_1(z)$.

Combined above assumptions with Eqs.(\ref{eqI1}), (\ref{eqI2}), (\ref{eqI3}) and (\ref{eqI4}), we know that $I_1=\emptyset$, $I_2=\{1, \cdots, l\}$, $I_3=\{l+1, \cdots, k\}$ and $I_4=\{k+1, \cdots, n-1\}$.

Let $I_2'=I_2\backslash\{1,l\}$ and $I_3'=I_3\backslash \{k\}$. It is seen that $I_2'=\emptyset$ when $l=1$.

\vskip 2mm
{\bf Subcase 3.1}. $j=\overline{n}$, i.e., $w\in V(G^{\overline{n}})$.

For $i\in [n-1]$, there exists a vertex $w_i\in V(G^{\overline{n}})$ whose out-neighbour $\widehat{w}_i\in V(G^i)$. By Lemma \ref{lemKfan}, there is an ($n-1$)-fan $P_1, \cdots, P_{n-1}$ in $G^{\overline{n}}$ from $w$ to $\{w_1, \cdots, w_{n-1}\}$ such that $w_i\in V(P_i)$ for $i\in [n-1]$.

\vskip 2mm
{\bf Subcase 3.1.1}. $l\ne 1$.

It is seen that there is a  $\{\widehat{\Gamma}_1(x), \widehat{w}_1, \widehat{\Gamma}_l(y), \widehat{\Gamma}_l(z)\}$-tree $\widehat{T}_1$ in $G^1 \cup G^{\overline{l}}$ according to Lemma \ref{lemBPn1}(5). Analogously, there is a $\{\widehat{x}, \widehat{y}, \widehat{w}_l\}$-tree $\widehat{T}_l$ in $G^{\overline{1}}\cup G^l$ and a $\{\widehat{\Gamma}_k(y), \widehat{z}, \widehat{w}_k\}$-tree $\widehat{T}_k$ in $G^{k}$, respectively.

\begin{figure}[htbp]
\begin{minipage}[t]{0.43\linewidth}
\centering
\resizebox{1\textwidth}{!} {\includegraphics{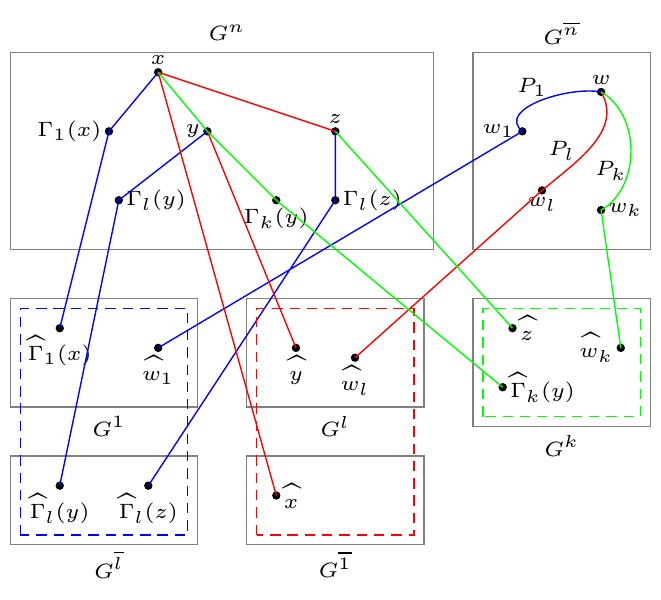}}
\caption{\small Illustration for Subcase 3.1.1} \label{figS3C31}
\end{minipage}
\begin{minipage}[t]{0.56\linewidth}
\centering
\resizebox{1\textwidth}{!} {\includegraphics{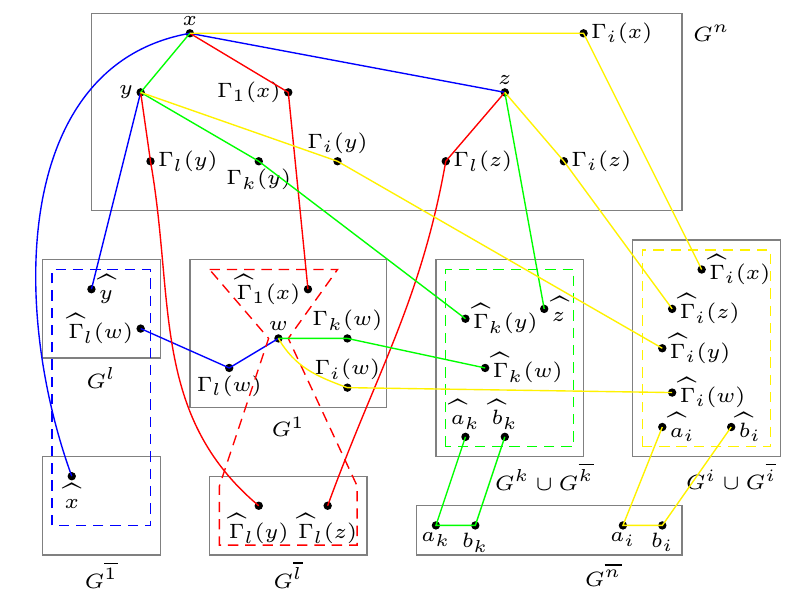}}
\caption{\small Illustration for Subcase 3.2.1.1} \label{figS3C32}
\end{minipage}
\end{figure}

Let
$$T_1=\widehat{T}_1\cup P(x, \widehat{\Gamma}_1(x)) \cup P(y, \widehat{\Gamma}_l(y))\cup P(z, \widehat{\Gamma}_l(z))\cup P_1\cup\{w_1\widehat{w}_1\},$$
$$T_l=\widehat{T}_l\cup P_l\cup\{x\widehat{x}, xz, y\widehat{y}, w_l\widehat{w}_l\},$$
and
$$T_k=\widehat{T}_k\cup P_k\cup P(y, \widehat{\Gamma}_k(y))\cup \{xy, z\widehat{z}, w_k\widehat{w}_k\}.$$
\noindent The trees $T_1$, $T_l$ and $T_k$ are depicted in Figure \ref{figS3C31} by blue, red and green lines, respectively.

There are vertices $a_i$ and $b_i=a_i(1)$ in $V(G^{\overline{k}})$ satisfying $\widehat{a}_i\in V(G^i)$ and $\widehat{b}_i\in V(G^{\overline{i}})$ for each $i\in I_2'\cup I_3'$. By similar arguments in Subcase 1.1.1, there is a $\{\widehat{\Gamma}_i(x), \widehat{\Gamma}_i(y), \widehat{\Gamma}_i(z), \widehat{w}_i\}$-inclusive tree $IT_i$ connects $G^i$ and $G^{\overline{i}}$ passing through $a_ib_i$, $i\in I_2'\cup I_3'$. For $i\in I_4$, there is a $\{\widehat{\Gamma}_i(x), \widehat{\Gamma}_i(y), \widehat{\Gamma}_i(z), \widehat{w}_i\}$-tree $\widehat{T}_i$ in $G^i$ since $G^i$ is connected.

Let $T_i$ be the same as in Eq.(\ref{eqTi6}) for $i\in I_2'\cup I_3'$ and $T_i$ be the same as in Eq.(\ref{eqTi5}) for $i\in I_4$. Then $T_1, \cdots,  T_{n-1}$ are ($n-1$) IDSTs in $BP_n$.

\vskip 2mm
{\bf Subcase 3.1.2}. $l=1$.

Remind that $I_2'=\emptyset$ under this situation. There are vertices $a_i$ and $b_i=a_i(1)$ in $V(G^{\overline{k}})$ satisfying $\widehat{a}_i\in V(G^i)$ and $\widehat{b}_i\in V(G^{\overline{i}})$ for each $1\le i\le k-1$. By analogous analysis in Subcase 3.1.1, there is a $\{\widehat{x}, \widehat{y}, \widehat{w}_1\}$-inclusive tree $IT_1$ connects $G^1$ and $G^{\overline{1}}$ passing through $a_1b_1$. Set
$$T_1=IT_1\cup P_1 \cup \{x\widehat{x}, y\widehat{y}, xz, w_1\widehat{w}_1\}.$$
For $i\in I_3\cup I_4$, let $T_i$ be the same as in Subcase 3.1.1. Then $T_1, \cdots,  T_{n-1}$ are ($n-1$) IDSTs in $BP_n$.

\vskip 2mm
{\bf Subcase 3.2}. $j\in \{1, \overline{1}, l, \overline{l}\}$.

We consider the case that $j=1$, for other cases the proof is similar. Since $w\in V(G^1)$, either $\widehat{\Gamma}_i(w)\in V(G^i)$ or $\widehat{\Gamma}_i(w)\in V(G^{\overline{i}})$ for $2\le i\le n-1$.

Set $W=\{\Gamma_2(w),\cdots, \Gamma_{n-1}(w)\}$ and $\widehat{W}=\{\widehat{\Gamma}_2(w),\cdots, \widehat{\Gamma}_{n-1}(w)\}.$

\vskip 2mm
{\bf Subcase 3.2.1}. $l\ne 1$.

\vskip 2mm
{\bf Subcase 3.2.1.1}. $\widehat{\Gamma}_l(w)\in V(G^l)$.

By analogous analysis in former discussions, $G^1\backslash W$ is connected. Note that $\widehat{W}\cap V(G^{\overline{l}})=\emptyset$. Thus, there is a $\{\widehat{\Gamma}_1(x), w, \widehat{\Gamma}_l(y), \widehat{\Gamma}_l(z)\}$-tree $\widehat{T}_1$ in $\big(G^1\backslash W\big)\cup G^{\overline{l}}$ since $\big(G^1\backslash W\big)\cup G^{\overline{l}}$ is connected by Lemma \ref{lemBPn1}(5). Furthermore, there is a $\{\widehat{x}, \widehat{y}, \widehat{\Gamma}_l(w)\}$-tree $\widehat{T}_l$ in $G^l\cup G^{\overline{1}}$.

Let
$$T_1=\widehat{T}_1\cup P(x, \widehat{\Gamma}_1(x))\cup P(y, \widehat{\Gamma}_l(y))\cup P(z, \widehat{\Gamma}_l(z))\}$$
and
$$T_l=\widehat{T}_l\cup  P(w, \widehat{\Gamma}_l(w))\cup \{x\widehat{x}, xz, y\widehat{y}\}.$$

For $i\in [n-1]\backslash \{1,l\}$, there are vertices $a_i$ and $b_i=a_i(1)$ in $V(G^{\overline{n}})$ satisfying $\widehat{a}_i\in V(G^i)$ and $\widehat{b}_i\in V(G^{\overline{i}})$. Due to similar analysis in Subcase 1.1.1, there is a $\{\widehat{z}, \widehat{\Gamma}_k(y), \widehat{\Gamma}_k(w)\}$-inclusive tree $IT_k$ connects $G^k$ and $G^{\overline{k}}$ passing through $a_kb_k$, furthermore, there is a $\{\widehat{\Gamma}_i(x), \widehat{\Gamma}_i(y), \widehat{\Gamma}_i(z), \widehat{\Gamma}_i(w)\}$-inclusive tree $IT_i$ connects $G^i$ and $G^{\overline{i}}$ passing through $a_ib_i$ for $i\in I_2'\cup I_3'\cup I_4$.

Let
$$T_k=IT_k \cup  P(y, \widehat{\Gamma}_k(y))\cup  P(w, \widehat{\Gamma}_k(w))\cup \{xy, z\widehat{z}\}.$$

For $i\in I_2'\cup I_3'\cup I_4$, let $T_i$ be the same as in Eq.(\ref{eqTi7}). See Figure \ref{figS3C32}. Then $T_1, \cdots,  T_{n-1}$ are ($n-1$) IDSTs in $BP_n$.

\vskip 2mm
{\bf Subcase 3.2.1.2}. $\widehat{\Gamma}_l(w)\in V(G^{\overline{l}})$.

There is a $\{\widehat{\Gamma}_1(x), w, \widehat{y}\}$-tree $\widehat{T}_1$ in $\big(G^1\backslash W\big)\cup G^{l}$ and a $\{\widehat{x}, \widehat{\Gamma}_l(y), \widehat{\Gamma}_l(z), \widehat{\Gamma}_l(w)\}$-tree $\widehat{T}_l$ in $G^{\overline{1}}\cup G^{\overline{l}}$.

Let
$$T_1=\widehat{T}_1\cup P(x, \widehat{\Gamma}_1(x))\cup \{y\widehat{y}, xz\}$$
and
$$T_l=\widehat{T}_l\cup P(y, \widehat{\Gamma}_l(y))\cup P(z, \widehat{\Gamma}_l(z))\cup P(w, \widehat{\Gamma}_l(w))\cup \{x\widehat{x}\}.$$

\begin{figure}[htbp]

\begin{minipage}[t]{0.5\linewidth}
\centering
\resizebox{1\textwidth}{!} {\includegraphics{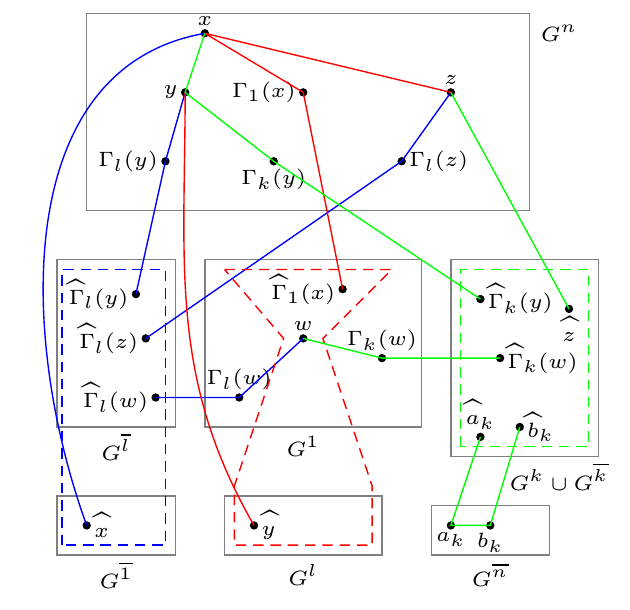}}
\caption{\small Illustration for Subcase 3.2.1.2} \label{figS3C33}
\end{minipage}
\begin{minipage}[t]{0.48\linewidth}
\centering
\resizebox{0.98\textwidth}{!} {\includegraphics{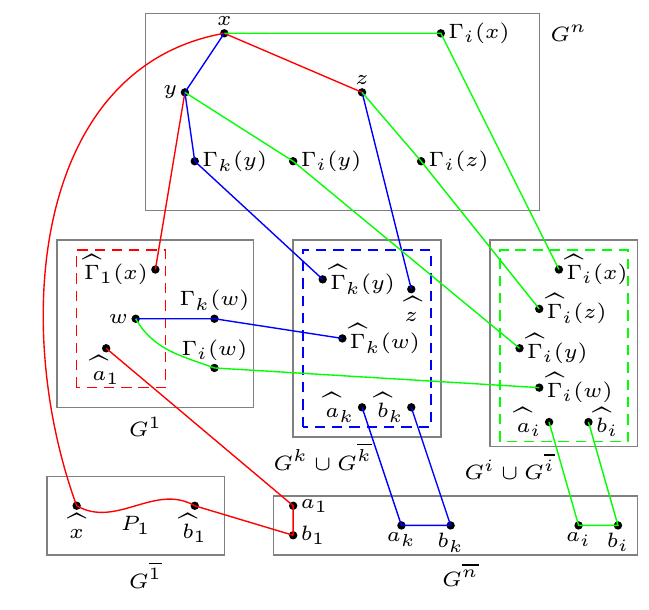}}
\caption{\small Illustration for Subcase 3.2.2} \label{figS3C322}
\end{minipage}
\end{figure}

\noindent See Figure \ref{figS3C33}. For $i\in [n-1]\backslash \{1,l\}$, let $T_i$ be the same as in Subcase 3.2.1.1.  Then $T_1, \cdots,  T_{n-1}$ are ($n-1$) IDSTs in $BP_n$.

\vskip 2mm
{\bf Subcase 3.2.2}. $l=1$.

There are vertices $a_i$ and $b_i=a_i(1)$ in $V(G^{\overline{n}})$ satisfying $\widehat{a}_i\in V(G^i)$ and $\widehat{b}_i\in V(G^{\overline{i}})$ for $i\in [n-1]$.
Remind that $G^1\backslash W$ is connected and $\widehat{\Gamma}_1(x)=\widehat{y}$, there is a $\{\widehat{y}, w, \widehat{a}_1\}$-tree $\widehat{T}_1$ in $G^1\backslash W$ and a ($\widehat{x}, \widehat{b}_1$)-path $P_1$ in $G^{\overline{1}}$. In addition, there is a $\{\widehat{\Gamma}_i(x), \widehat{\Gamma}_i(y), \widehat{\Gamma}_i(z), \widehat{\Gamma}_i(w)\}$-inclusive tree $IT_i$ connects $G^i$ and $G^{\overline{i}}$ passing through $a_ib_i$ for $2\le i\le n-1$.

Let $T_1=\widehat{T}_1\cup P_1\cup \{x\widehat{x}, y\widehat{y}, xz, a_1\widehat{a}_1, a_1b_1, b_1\widehat{b}_1\}$ and let $T_i$ be the same as in Subcase 3.2.1.1 for $i\in I_3\cup I_4$. See Figure \ref{figS3C322}. Then $T_1, \cdots,  T_{n-1}$ are ($n-1$) IDSTs in $BP_n$.

\vskip 2mm
{\bf Subcase 3.3}. $j=k$, i.e., $w\in V(G^k)$.

Let $W=\{\Gamma_i(w)| i\in I_2\cup I_3'\cup I_4\}$.

Clearly, $\{\widehat{z}, \widehat{\Gamma}_k(y)\}\cap W=\emptyset.$ There is a $\{\widehat{z}, \widehat{\Gamma}_k(y), w\}$-tree $\widehat{T}_k$ in $G^k\backslash W$ since $G^k\backslash W$ is connected. Let $T_k=\widehat{T}_k\cup P(y, \widehat{\Gamma}_k(y))\cup\{z\widehat{z}, xy\}$.

\vskip 2mm
{\bf Subcase 3.3.1}. $l\ne 1$.

For $i\in I_2'\cup I_3'\cup I_4$, there exist vertices $a_i$ and $b_i=a_i(1)$ in $V(G^{\overline{n}})$ satisfying $\widehat{a}_i\in V(G^i)$ and $\widehat{b}_i\in V(G^{\overline{i}})$.  For $i\in I_2'\cup I_3'\cup I_4$, we can find a $\{\widehat{\Gamma}_i(x), \widehat{\Gamma}_i(y), \widehat{\Gamma}_i(z), \widehat{\Gamma}_i(w)\}$-inclusive tree $IT_i$ connects $G^i$ and $G^{\overline{i}}$ passing through $a_ib_i$ and let $T_i$ be the same as in Eq.(\ref{eqTi7}).

Consider the locations of $\widehat{\Gamma}_1(w)$ and $\widehat{\Gamma}_l(w)$. To avoid duplication, we consider the following two cases.

\vskip 2mm
{\bf Subcase 3.3.1.1}. $\widehat{\Gamma}_1(w)\in V(G^1)$ and $\widehat{\Gamma}_l(w)\in V(G^l)$.

There is a $\{\widehat{\Gamma}_1(x), \widehat{\Gamma}_1(w), \widehat{\Gamma}_l(y), \widehat{\Gamma}_l(z)\}$-tree $\widehat{T}_1$ in $G^1\cup G^{\overline{l}}$ and a $\{\widehat{x}, \widehat{y}, \widehat{\Gamma}_l(w)\}$-tree  $\widehat{T}_l$ in $G^l\cup G^{\overline{1}}$, respectively.

\begin{figure}[htbp]
\begin{minipage}[t]{0.49\linewidth}
\centering
\resizebox{0.9\textwidth}{!} {\includegraphics{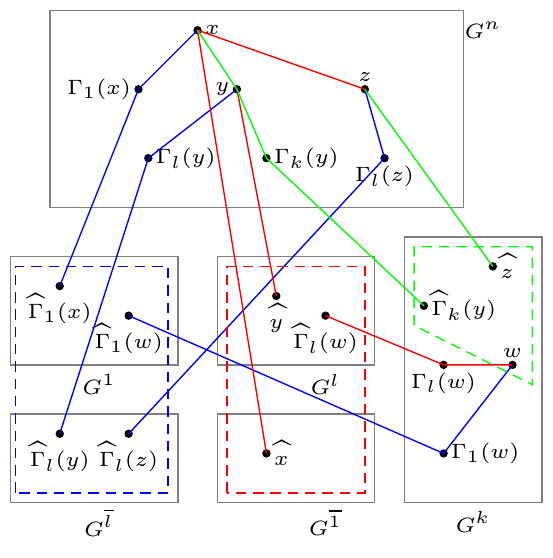}}
\caption{\small Illustration for Subcase 3.3.1.1} \label{figS3C34}
\end{minipage}
\begin{minipage}[t]{0.49\linewidth}
\centering
\resizebox{0.9\textwidth}{!} {\includegraphics{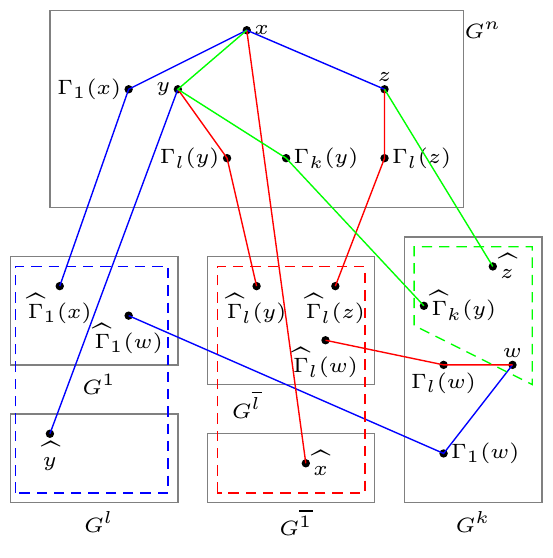}}
\caption{\small Illustration for Subcase 3.3.1.2} \label{figS3C35}
\end{minipage}
\end{figure}

Let
$$T_1=\widehat{T}_1\cup P(x, \widehat{\Gamma}_1(x))\cup P(y, \widehat{\Gamma}_l(y)) \cup P(z, \widehat{\Gamma}_l(z))\cup P(w, \widehat{\Gamma}_1(w))$$
and
$$T_l=\widehat{T}_l\cup  P(w, \widehat{\Gamma}_l(w))\cup \{x\widehat{x}, y\widehat{y}, xz\}.$$
For simplicity, only the trees $T_1$, $T_l$ and $T_k$ are depicted in Figure \ref{figS3C34}. Then $T_1, \cdots,  T_{n-1}$ are ($n-1$) IDSTs in $BP_n$.

\vskip 2mm
{\bf Subcase 3.3.1.2}. $\widehat{\Gamma}_1(w)\in V(G^1)$ and $\widehat{\Gamma}_l(w)\in V(G^{\overline{l}})$.

There is a $\{\widehat{\Gamma}_1(x), \widehat{\Gamma}_1(w), \widehat{y}\}$-tree $\widehat{T}_1$ in $G^1\cup G^{l}$ and a $\{\widehat{x}, \widehat{\Gamma}_l(y), \widehat{\Gamma}_l(z), \widehat{\Gamma}_l(w)\}$-tree  $\widehat{T}_l$ in $G^{\overline{l}}\cup G^{\overline{1}}$, respectively.

Let
$$T_1=\widehat{T}_1\cup P(x, \widehat{\Gamma}_1(x))\cup P(w, \widehat{\Gamma}_1(w))\cup \{y\widehat{y}, xz\}$$
and
$$T_l=\widehat{T}_l\cup  P(y, \widehat{\Gamma}_l(y))\cup P(z, \widehat{\Gamma}_l(z))\cup P(w, \widehat{\Gamma}_l(w))\cup \{x\widehat{x}\}.$$
See Figure \ref{figS3C35}. Then $T_1, \cdots,  T_{n-1}$ are ($n-1$) IDSTs in $BP_n$.

\vskip 2mm
{\bf Subcase 3.3.2}. $l=1$.

For $i\in I_2\cup I_3'\cup I_4$, there exist vertices $a_i$ and $b_i=a_i(1)$ in $V(G^{\overline{n}})$ satisfying $\widehat{a}_i\in V(G^i)$ and $\widehat{b}_i\in V(G^{\overline{i}})$. Hence, there is a $\{\widehat{x}, \widehat{y}, \widehat{\Gamma}_1(w)\}$-inclusive tree $IT_1$ connects $G^1$ and $G^{\overline{1}}$ passing through $a_1b_1$.

Let
$$T_1=IT_1\cup P(w, \widehat{\Gamma}_1(w))\cup \{x\widehat{x}, y\widehat{y}, xz\}$$
and let $T_i$ be the same as in Subcase 3.3.1 for $i\in I_3\cup I_4$. Then $T_1, \cdots,  T_{n-1}$ are ($n-1$) IDSTs in $BP_n$.

\vskip 2mm
{\bf Subcase 3.4}. $j=\overline{k}$.

Let $W=\{\Gamma_i(w)~|~ i\in [n-1]\backslash \{k\}\}$.

For $i\in [n-1]$, there are vertices $a_i$ and $b_i=a_i(1)$ in $V(G^{\overline{n}})$ satisfying $\widehat{a}_i\in V(G^i)$ and $\widehat{b}_i\in V(G^{\overline{i}})$.

Since $G^k\backslash W$ is connected, there exists a $\{\widehat{z}, \widehat{\Gamma}_k(y),w\}$-inclusive tree $IT_k$ connects $(G^k\backslash W)$ and $G^{\overline{k}}$ passing through $a_kb_k$.  Let
$$T_k=IT_k\cup P(y, \widehat{\Gamma}_k(y))\cup \{xy, z\widehat{z}\}.$$

For $i\in I_2\cup I_3'\cup I_4$, we can construct $T_i$ be the same as in Subcase 3.3. Then $T_1, \cdots,  T_{n-1}$ are ($n-1$) IDSTs in $BP_n$.

\vskip 2mm
{\bf Subcase 3.5}. Either $j\in I_2'\cup I_3'\cup I_4$ or $\overline{j}\in I_2'\cup I_3'\cup I_4$.

We may construct $T_1$ and $T_l$ by analogous arguments in Subcase 3.3 and let $T_i$ ($i\in [n-1]\backslash \{1,l\}$) by similar methods in Subcase 3.2. Then ($n-1$) IDSTs in $BP_n$ can be obtained. The proof is completed. \hfill$\Box$

\section{$\max\{|S\cap V(G^i)|\}=2$ for $i\in [[n]]$}\label{secS2}

\begin{lemma}\label{lemK4S22}
For $n\ge 3$, let $BP_n=G^1\oplus G^{\overline{1}}\oplus G^2 \oplus G^{\overline{2}}\oplus \cdots \oplus  G^n \oplus G^{\overline{n}}$ and $S=\{x,y,z,w\}$ be any 4-subset of $V(BP_n)$. If there are different integers $i$ and $j$ in $[[n]]$ that $|S\cap V(G^i)|=|S\cap V(G^j)|=2$. Then there exist ($n-1$) IDSTs in $BP_n$.
\end{lemma}

\noindent{\bf Proof}\;  Without loss of generality, we may assume that $\{x,y\}\subseteq V(G^n)$ and $\{z,w\}\subseteq V(G^j)$, where $j\in [[n]]\backslash \{n\}$. Furthermore, we may assume that $x=12\cdots n$.

By Lemma \ref{lemxypath} and the fact that $\kappa(G^n)=n-1$, there is a family of ($n-1$) IDPs $P_1, \cdots, P_{n-1}$ in $G^n$ connecting $x$ and $y$.  Without loss of generality, assume that $x(i)\in V(P_i)$ for $i\in [n-1]$. It is possible that $y=x(i)$ for $i\in [n-1]$. This possibility does not affect the following discussions. Note that $\widehat{x}(i)\in V(G^i)$ for $i\in [n-1]$ and $\widehat{x}\in V(G^{\overline{1}})$.

Similarly, there is a family of ($n-1$) IDPs $Q_1, \cdots, Q_{n-1}$ in $G^j$ connecting $z$ and $w$. For simplicity, denote $q_i(z)$ be the neighbour of $z$ which belongs to $V(Q_i)$ for $i\in [n-1]$.

Let $\mathcal{Q}=\{\widehat{q}_i(z)~|~ i\in [n-1]\}$ and $\mathcal{Q}'=\mathcal{Q}\cup \{\widehat{z}\}$.

By Lemma \ref{lemBPn1}(4), $|\mathcal{Q}'\cap V(G^i)|\le 1$ for each $i\in [[n]]$. Especially, $|\mathcal{Q}'\cap V(G^n)|\le 1$. Without loss of generality, we may assume that $\mathcal{Q}\cap V(G^n)=\emptyset$. We may use $\widehat{z}$ to replace $\widehat{q}_i(z)$ in the following discussions if $\widehat{q}_i(z)\in V(G^n)$ for $i\in [n-1]$.

\vskip 2mm
{\bf Case 1}. $j\notin [n-1]$.

Let $I_s=\{i~|~ i\in [n-1] ~{\rm and}~ |\mathcal{Q}\cap V(G^i)|= 1\}.$

Without loss of generality, we may assume that $I_s=\{1,\cdots, k\}$ and $\widehat{q}_i(z)\in V(G^i)$ for $1\le i\le k$. When $k<n-2$, for $k+1\le i\le n-1$, we may relabel the vertices $q_i(z)$ that $\widehat{q}_i(z)\in V(G^{m_i})$, where $m_i$ is an integer different from $\overline{i}$.

For $1\le i\le k$, there is a ($\widehat{x}(i),\widehat{q}_i(z)$)-path $R_i$ in $G^i$ since $G^i$ is connected. Moreover, for $k+1\le i\le n-1$, there is a ($\widehat{x}(i),\widehat{q}_i(z)$)-path $R_i$ in $G^i\cup G^{m_i}$ by Lemma \ref{lemBPn1}(5) when $k<n-2$. If $k=n-2$. Set $G'=G^n\cup G^j\cup \big(\bigcup_{i=1}^{n-2}G^i \big).$  Combined Lemma \ref{lemBPn1}(5) with the fact that there are $n\ge 3$ clusters in $BP_n\backslash G'$, $BP_n\backslash G'$ is connected. Note that $\{\widehat{x}(n-1),\widehat{q}_{n-1}(z)\}\subseteq V(G')$. Thus, there is a ($\widehat{x}(n-1),\widehat{q}_{n-1}(z)$)-path $R_{n-1}$ in $BP_n\backslash G'$. Let
$$T_i=P_i\cup Q_i\cup R_i\cup \{x(i)\widehat{x}(i), q_i(z)\widehat{q}_i(z)\}, \;\;\; 1\le i\le n-1.$$
Then $T_1, \cdots,  T_{n-1}$ are ($n-1$) IDSTs in $BP_n$.

\vskip 2mm
{\bf Case 2}. $j\in [n-1]$.

To avoid duplication, we may assume that $j=n-1$.

Since $G^{n-1}$ is connected, there is a path $P$ in $G^{n-1}$ connecting $\widehat{x}(n-1)$ and $z$. Let $v$ be the first vertex in $V(P)\cap \Big(\bigcup_{i=1}^{n-1}V(Q_i)\Big)$ when $P$ starts from $\widehat{x}(n-1)$. Without loss of generality, assume that $v\in V(Q_{n-1})$. Denote by $P(\widehat{x}(n-1),v)$ the subgraph of $P$ that starts from $\widehat{x}(n-1)$ and ends at $v$.

Let
$$T_{n-1}=P_{n-1}\cup Q_{n-1}\cup P(\widehat{x}(n-1),v)\cup \{x(n-1)\widehat{x}(n-1)\}$$
and construct $T_1, \cdots, T_{n-2}$ by similar methods in Case 1. Then $T_1, \cdots,  T_{n-1}$ are ($n-1$) IDSTs in $BP_n$. \hfill$\Box$

\begin{lemma}\label{lemK4S211}
For $n\ge 3$, let $BP_n=G^1\oplus G^{\overline{1}}\oplus G^2 \oplus G^{\overline{2}}\oplus \cdots \oplus  G^n \oplus G^{\overline{n}}$ and $S=\{x,y,z,w\}$ be any 4-subset of $V(BP_n)$. If there are different integers $i,j$ and $l$ in $[[n]]$ such that $|S\cap V(G^i)|=2$ and $|S\cap V(G^j)|=|S\cap V(G^l)|=1$. Then there exist ($n-1$) IDSTs in $BP_n$.
\end{lemma}

\noindent
{\bf Proof} \; Without loss of generality, we may assume that $\{x,y\}\subseteq V(G^n)$, $z\in V(G^j)$ and $w\in V(G^l)$, where $j$ and $l$ are different integers in $[[n]]\backslash\{n\}$. Furthermore, let $x=12\cdots n$.

Since $\kappa(G^n)=n-1$, there is a family of ($n-1$) IDPs $P_1, \cdots, P_{n-1}$ in $G^n$ connecting $x$ and $y$. For simplicity, assume that $x(i)\in V(P_i)$ for $i\in [n-1]$. Note that $\widehat{x}\in V(G^{\overline{1}})$ and $\widehat{x}(i)\in V(G^{i})$ for $i\in [n-1]$.

The following discussions are based on the values of $j$ and $l$.

\vskip 2mm
{\bf Case 1}. $\{j,l\}\cap [n-1]=\emptyset$.

That means $\{j,l\}\subseteq \{\overline{1}, \overline{2}, \cdots, \overline{n}\}.$

\vskip 2mm
{\bf Subcase 1.1}. $j=\overline{n}$ or $l=\overline{n}$.

Without loss of generality, assume that $z\in V(G^{\overline{n}})$ and $w\in V(G^{\overline{n-1}})$.

According to Lemma \ref{lemBPn1}(5), there are ($n-1$) different vertices $u_1, \cdots, u_{n-1}$ in $V(G^{\overline{n}})$ that $\widehat{u}_i\in V(G^{i})$ for $i\in [n-1]$. By Lemma \ref{lemKfan}, there is an ($n-1$)-fan $Q_1, \cdots, Q_{n-1}$ in $G^{\overline{n}}$ from $z$ to $\{u_1, \cdots, u_{n-1}\}$ such that $u_i\in V(Q_i)$ for $i\in [n-1]$.

Likewise, there are ($n-1$) different vertices $v_1, \cdots, v_{n-1}$ in $V(G^{\overline{n-1}})$ that $\widehat{v}_i\in V(G^{i})$ for $i\in [n-2]$ and $\widehat{v}_{n-1}\in V(G^{\overline{1}})$. Moreover, there is an ($n-1$)-fan $R_1, \cdots, R_{n-1}$  in $G^{\overline{n-1}}$ from $w$ to $\{v_1, \cdots, v_{n-1}\}$ such that $v_i\in V(R_i)$ for $i\in [n-1]$.

By Lemma \ref{lemBPn1}, there is a $\{\widehat{x}(n-1), \widehat{u}_{n-1}, \widehat{v}_{n-1}\}$-tree $\widehat{T}_{n-1}$ in $G^{n-1}\cup G^{\overline{1}}$ and a $\{\widehat{x}(i), \widehat{u}_i, \widehat{v}_i\}$-tree $\widehat{T}_i$ in $G^i$, $i\in [n-2]$.

\begin{figure}[htbp]
\begin{minipage}[t]{0.47\linewidth}
\centering
\resizebox{0.8\textwidth}{!} {\includegraphics{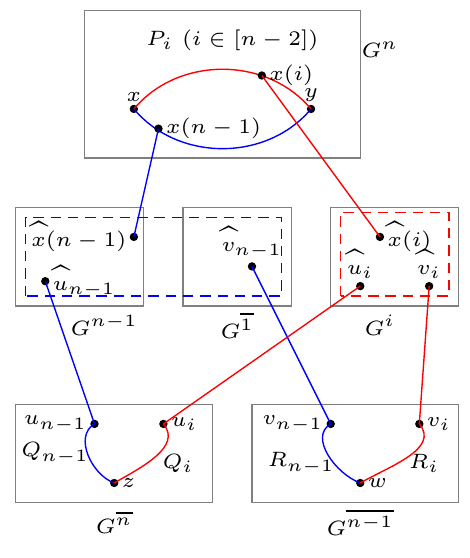}}
\caption{\small Illustration for Subcase 1.1} \label{figS211C11}
\end{minipage}
\begin{minipage}[t]{0.5\linewidth}
\centering
\resizebox{0.9\textwidth}{!} {\includegraphics{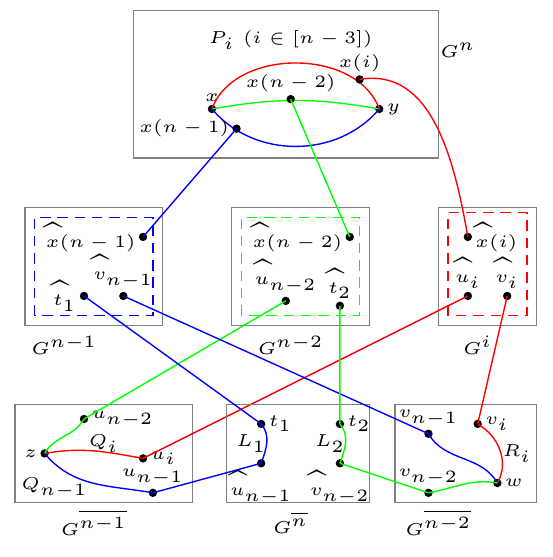}}
\caption{\small Illustration for Subcase 1.2} \label{figS211C12}
\end{minipage}
\end{figure}

For $i\in [n-1]$, let
\begin{eqnarray}\label{EqTi8}
T_i=\widehat{T}_i\cup P_i\cup Q_i\cup R_i\cup \{x(i)\widehat{x}(i), u_i\widehat{u}_i, v_i\widehat{v}_i\}.
\end{eqnarray}
See Figure \ref{figS211C11}. Then $T_1, \cdots,  T_{n-1}$ are ($n-1$) IDSTs in $BP_n$.

\vskip 2mm
{\bf Subcase 1.2}. $j\ne\overline{n}$ and $l\ne\overline{n}$.

Without loss of generality, assume that $z\in V(G^{\overline{n-1}})$ and $w\in V(G^{\overline{n-2}})$.

By Lemma \ref{lemBPn1}, there are ($n-1$) different vertices $u_1, \cdots, u_{n-1}$ in $V(G^{\overline{n-1}})$ that $\widehat{u}_i\in V(G^{i})$ for $i\in [n-2]$ and $\widehat{u}_{n-1}\in V(G^{\overline{n}})$. By Lemma \ref{lemKfan}, there is an ($n-1$)-fan $Q_1, \cdots, Q_{n-1}$  in $G^{\overline{n-1}}$ from $z$ to $\{u_1, \cdots, u_{n-1}\}$ such that $u_i\in V(Q_i)$ for $i\in [n-1]$.

Analogously, there are ($n-1$) different vertices $v_1, \cdots, v_{n-1}$ in $V(G^{\overline{n-2}})$ that $\widehat{v}_i\in V(G^{i})$ for $ i\in [n-1]\backslash\{n-2\}$ and $\widehat{v}_{n-2}\in V(G^{\overline{n}})$. Moreover, there is an ($n-1$)-fan $R_1, \cdots, R_{n-1}$  in $G^{\overline{n-2}}$ from $w$ to $\{v_1, \cdots, v_{n-1}\}$ such that $v_i\in V(R_i)$ for $i\in [n-1]$.

By Lemma \ref{lemBPn2}, let $t_1$ be the neighbour of $\widehat{u}_{n-1}$ in $V(G^{\overline{n}})$ such that $\widehat{t}_1\in V(G^{n-1})$. Likewise, let $t_2$ be the neighbour of $\widehat{v}_{n-2}$ in $V(G^{\overline{n}})$ that $\widehat{t}_2\in V(G^{n-2})$. For simplicity, denote by $L_1$ the path $\{u_{n-1}\widehat{u}_{n-1}, t_{1}\widehat{u}_{n-1}, t_1\widehat{t}_1\}$ and denote by $L_2$ the path $\{v_{n-2}\widehat{v}_{n-2}, t_{2}\widehat{v}_{n-2}, t_2\widehat{t}_2\}$.

Note that there is a  $\{\widehat{x}(n-2), \widehat{u}_{n-2}, \widehat{t}_2\}$-tree $\widehat{T}_{n-2}$ in $G^{n-2}$ and a $\{\widehat{x}(n-1), \widehat{v}_{n-1}, \widehat{t}_1\}$-tree $\widehat{T}_{n-1}$ in $G^{n-1}$, respectively.

Let
$$T_{n-2}=\widehat{T}_{n-2}\cup P_{n-2}\cup Q_{n-2}\cup R_{n-2}\cup L_2\cup\{x(n-2)\widehat{x}(n-2), u_{n-2}\widehat{u}_{n-2} \}$$
and
$$T_{n-1}=\widehat{T}_{n-1}\cup P_{n-1}\cup Q_{n-1}\cup R_{n-1}\cup L_1\cup\{x(n-1)\widehat{x}(n-1),  v_{n-1}\widehat{v}_{n-1}\}.$$
For $1\le i\le n-3$, let $\widehat{T}_i$ and $T_i$ be the same as in Subcase 1.1. See Figure \ref{figS211C12}. Then $T_1, \cdots,  T_{n-1}$ are ($n-1$) IDSTs in $BP_n$.

\vskip 2mm
{\bf Case 2}. $|\{j,l\}\cap [n-1]|=1$.

Without loss of generality, we may assume that $j=n-1$. That means $z\in V(G^{n-1})$.

There exist ($n-2$) different vertices $u_1, \cdots, u_{n-2}$ in $V(G^{n-1})\backslash \{\widehat{x}(n-1)\}$ such that $\widehat{u}_i\in V(G^{i})$ for $i\in [n-2]$. By Lemma \ref{lemKfan}, there is an ($n-1$)-fan $Q_1, \cdots, Q_{n-1}$ in $G^{n-1}$ from $z$ to $\{u_1, \cdots, u_{n-2}, \widehat{x}(n-1)\}$ such that $u_i\in V(Q_i)$ for $i\in [n-2]$ and $\widehat{x}(n-1)\in V(Q_{n-1})$.

\vskip 2mm
{\bf Subcase 2.1}. $l=\overline{1}$, i.e., $w\in V(G^{\overline{1}})$.

There exist ($n-2$) different vertices $v_1, \cdots, v_{n-2}$ in $V(G^{\overline{1}})\backslash \{\widehat{x}\}$ such that $\widehat{v}_i\in V(G^{i})$ for $2\le i\le n-2$ and $\widehat{v}_1\in V(G^{\overline{n}})$. Furthermore, there is an ($n-1$)-fan $R_1, \cdots, R_{n-1}$ in $G^{\overline{1}}$ from $w$ to $\{v_1, \cdots, v_{n-2}, \widehat{x}\}$ such that $v_i\in V(R_i)$ for $i\in [n-2]$ and $\widehat{x}\in V(R_{n-1})$.

Let $T_{n-1}=P_{n-1}\cup Q_{n-1}\cup R_{n-1}\cup \{x\widehat{x}, x(n-1)\widehat{x}(n-1)\}$.

\begin{figure}[htbp]
\begin{minipage}[t]{0.49\linewidth}
\centering
\resizebox{0.9\textwidth}{!} {\includegraphics{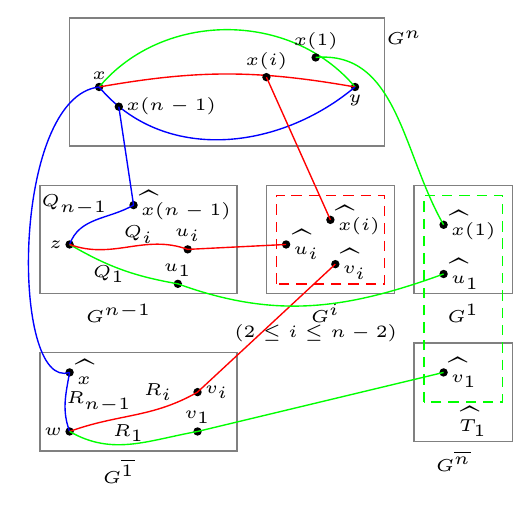}}
\caption{\small Illustration for Subcase 2.1} \label{figS211C21}
\end{minipage}
\begin{minipage}[t]{0.49\linewidth}
\centering
\resizebox{0.9\textwidth}{!} {\includegraphics{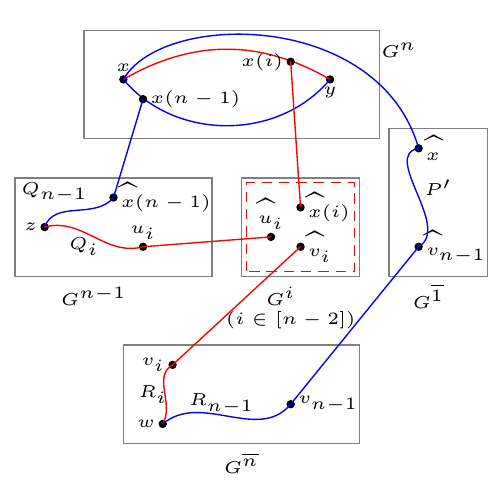}}
\caption{\small Illustration for Subcase 2.2} \label{figS211C22}
\end{minipage}
\end{figure}

Note that there is a $\{\widehat{x}(1), \widehat{u}_1, \widehat{v}_1\}$-tree $\widehat{T}_1$ in $G^1\cup G^{\overline{n}}$ and a $\{\widehat{x}(i), \widehat{u}_i, \widehat{v}_i\}$-tree $\widehat{T}_i$ in $G^i$ for $2\le i\le n-2$. We may construct $T_1,\cdots, T_{n-2}$ be the same as in Eq.(\ref{EqTi8}). See Figure \ref{figS211C21}. Then ($n-1$) IDSTs $T_1,\cdots, T_{n-1}$ can be obtained.

\vskip 2mm
{\bf Subcase 2.2}. $l\ne\overline{1}$.

To avoid duplication, we consider the case that $l=\overline{n}$. It follows that $w\in V(G^{\overline{n}})$.

There are ($n-1$) different vertices $v_1,\cdots, v_{n-1}$ in $V(G^{\overline{n}})$ such that $\widehat{v}_i\in V(G^i)$ for $i\in [n-2]$ and $\widehat{v}_{n-1}\in V(G^{\overline{1}})$. Furthermore, there is an ($n-1$)-fan $R_1, \cdots, R_{n-1}$ in $G^{\overline{n}}$ from $w$ to $\{v_1, \cdots, v_{n-1}\}$ such that $v_i\in V(R_i)$ for $i\in [n-1]$. Moreover, there is a ($\widehat{x},\widehat{v}_{n-1}$)-path $P'$ in $G^{\overline{1}}$ since $G^{\overline{1}}$ is connected.

Let $T_{n-1}=P_{n-1}\cup Q_{n-1}\cup R_{n-1}\cup P'\cup \{x\widehat{x}, x(n-1)\widehat{x}(n-1), v_{n-1}\widehat{v}_{n-1}\}$ and $T_i$ be the same as in Eq.(\ref{EqTi8}) for $i\in[n-2]$. See Figure \ref{figS211C22}. Then $T_1,\cdots, T_{n-1}$ are ($n-1$) IDSTs in $BP_n$.

\vskip 2mm
{\bf Case 3}. $|\{j,l\}\cap [n-1]|=2$, i.e., $\{j,l\}\subseteq [n-1]$.

Without loss of generality, we may assume that $j=1$ and $l=2$.

There are ($n-2$) vertices $u_i$ in $V(G^1)\backslash \{\widehat{x}(1)\}$ such that $\widehat{u}_i\in V(G^i)$ for $2\le i\le n-1$. By Lemma \ref{lemKfan}, there is an ($n-1$)-fan $Q_1, \cdots, Q_{n-1}$ in $G^1$ from $z$ to $\{\widehat{x}(1), u_2, \cdots, u_{n-1}\}$ such that $\widehat{x}(1)\in V(Q_1)$ and $u_i\in V(Q_i)$ for $2\le i\le n-1$.

Set $G'=BP_n\backslash \big(\bigcup_{i=1}^{n}V(G^i)\big)$. In order to obtain ($n-1$) IDSTs in $BP_n$, the following discussions are made according to the location of $\widehat{w}$.

\vskip 2mm
{\bf Subcase 3.1}. $\widehat{w}\in V(G')$.

There exist ($n-3$) different vertices $v_3,\cdots, v_{n-1}$ in $V(G^2)\backslash \{\widehat{x}(2),\widehat{u}_2\}$ such that $\widehat{v}_i\in V(G^i)$ for $3\le i\le n-1$. By Lemma \ref{lemKfan}, there is an ($n-1$)-fan $R_1, \cdots, R_{n-1}$ in $G^2$ from $w$ to $\{\widehat{u}_2, \widehat{x}(2), v_3,\cdots, v_{n-1}\}$ such that $\widehat{u}_2\in V(R_1)$, $\widehat{x}(2)\in V(R_2)$ and $v_i\in V(R_i)$ for $3\le i\le n-1$.

For $3\le i\le n-1$, there is a $\{\widehat{x}(i), \widehat{v}_i, \widehat{u}_i\}$-tree $\widehat{T}_i$ in $G^i$ since $\{\widehat{x}(i), \widehat{v}_i, \widehat{u}_i\}\subseteq V(G^i)$, thus, let $T_i$ be the same as in Eq.(\ref{EqTi8}).

\begin{figure}[htbp]
\begin{minipage}[t]{0.49\linewidth}
\centering
\resizebox{0.9\textwidth}{!} {\includegraphics{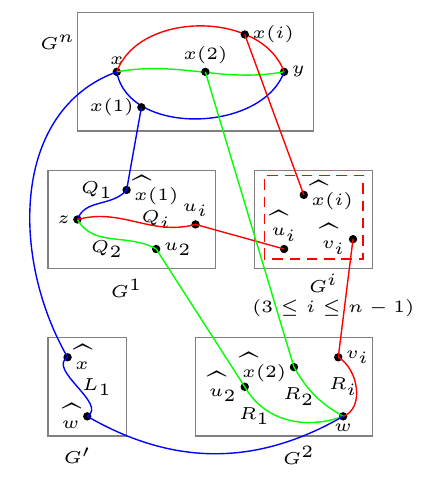}}
\caption{\small Illustration for Subcase 3.1} \label{figS211C31}
\end{minipage}
\begin{minipage}[t]{0.49\linewidth}
\centering
\resizebox{0.9\textwidth}{!} {\includegraphics{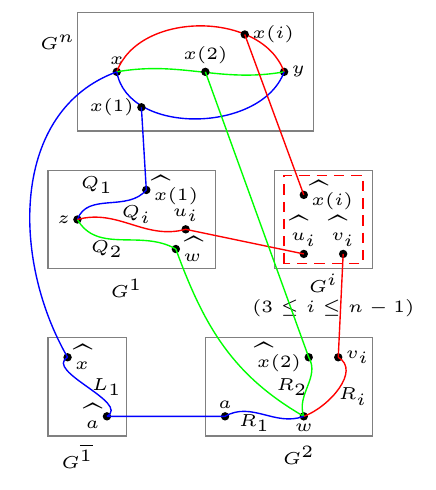}}
\caption{\small Illustration for Subcase 3.2} \label{figS211C32}
\end{minipage}
\end{figure}

It is seen that there is a ($\widehat{x},\widehat{w}$)-path $L_1$ in $G'$ by Lemma \ref{lemBPn1}(5). Let $T_2=P_2\cup R_2\cup R_1\cup Q_2\cup \{x(2)\widehat{x}(2), u_2\widehat{u}_2\}$ and $T_1=P_1\cup L_1\cup Q_1\cup \{x\widehat{x}, x(1)\widehat{x}(1), w\widehat{w}\}$. See Figure \ref{figS211C31}. Then $T_1,\cdots, T_{n-1}$ are ($n-1$) IDSTs in $BP_n$.

\vskip 2mm
{\bf Subcase 3.2}. $\widehat{w}\in V(G^i)$, for $i\in [n-1]\backslash\{2\}$.

Without loss of generality, we may assume that $\widehat{w}\in V(G^1)$. Let $u_2=\widehat{w}$ in former discussions. Then $\widehat{w}\in V(Q_2)$.

We can find ($n-2$) vertices $a, v_3,\cdots, v_{n-1}$ in $V(G^2)\backslash \{\widehat{x}(2)\}$ such that $\widehat{a}\in V(G^{\overline{1}})$ and $\widehat{v}_i\in V(G^i)$ for $3\le i\le n-1$. By Lemma \ref{lemKfan}, there is an ($n-1$)-fan $R_1, \cdots, R_{n-1}$ in $G^2$ from $w$ to $\{ a, \widehat{x}(2), v_3, \cdots, v_{n-1}\}$ that $a\in V(R_1)$, $\widehat{x}(2)\in V(R_2)$ and $v_i\in V(R_i)$ for $3\le i\le n-1$. Moreover, there is a ($\widehat{x},\widehat{a}$)-path $L_1$ in $G^{\overline{1}}$.

Let $T_2=P_2\cup R_2\cup Q_2\cup \{x(2)\widehat{x}(2), w\widehat{w}\}$, $T_1=P_1\cup L_1\cup Q_1\cup R_1\cup \{x\widehat{x}, x(1)\widehat{x}(1), a\widehat{a}\}$ and $T_i$ ($3\le i\le n-1$) be the same as in Eq.(\ref{EqTi8}). See Figure \ref{figS211C32}. Then $T_1,\cdots, T_{n-1}$ are ($n-1$) IDSTs in $BP_n$.

\vskip 2mm
{\bf Subcase 3.3}. $\widehat{w}\in V(G^n)$.

Since $G^n$ is connected, there is a path $\widetilde{P}$ in $G^n$ connecting $\widehat{w}$ and $x$. Let $v$ be the first vertex in $V(\widetilde{P})\cap \big(\bigcup_{i=1}^{n-1}V(P_i)\big)$ when $\widetilde{P}$ starts from $\widehat{w}$. Denote by $\widetilde{P}(\widehat{w},v)$ the subgraph of $\widetilde{P}$ that starts from $\widehat{w}$ and ends at $v$.

\vskip 2mm
{\bf Subcase 3.3.1}. $v\in V(P_1)$.

Let $T_1=P_1\cup Q_1\cup \widetilde{P}(\widehat{w},v)\cup \{w\widehat{w}, x(1)\widehat{x}(1)\}$ and $T_i$ be the same as in Subcase 3.1 for $2\le i\le n-1$. Then $T_1,\cdots, T_{n-1}$ are ($n-1$) IDSTs in $BP_n$.

\vskip 2mm
{\bf Subcase 3.3.2}. $v\in V(P_2)$.

There are ($n-2$) vertices $a, v_3, \cdots, v_{n-1}$ in $V(G^2)\backslash \{\widehat{x}(2), \widehat{u}_2\}$ such that $\widehat{a}\in V(G^{\overline{1}})$ and $\widehat{v}_i\in V(G^i)$ for $3\le i\le n-1$. Furthermore, there is an ($n-1$)-fan $R_1, \cdots, R_{n-1}$ in $G^2$ from $w$ to $\{\widehat{u}_2, a, v_3, \cdots, v_{n-1}\}$ that $\widehat{u}_2\in V(R_1)$, $a\in V(R_2)$ and $v_i\in V(R_i)$ for $3\le i \le n-1$. Moreover, there is a ($\widehat{x},\widehat{a}$)-path $L_1$ in $G^{\overline{1}}$.

\begin{figure}[htbp]

\begin{minipage}[t]{0.47\linewidth}
\centering
\resizebox{0.9\textwidth}{!} {\includegraphics{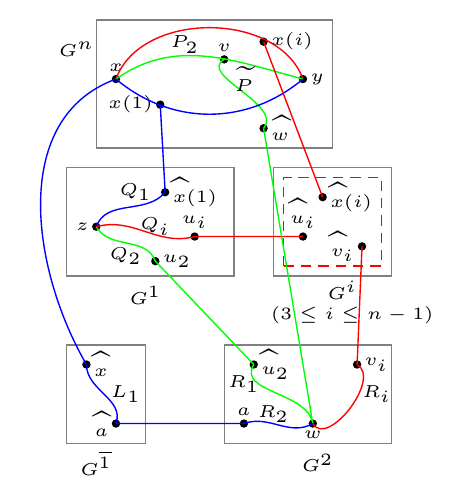}}
\caption{\small Illustration for Subcase 3.3.2} \label{figS211C332}
\end{minipage}
\begin{minipage}[t]{0.5\linewidth}
\centering
\resizebox{0.9\textwidth}{!} {\includegraphics{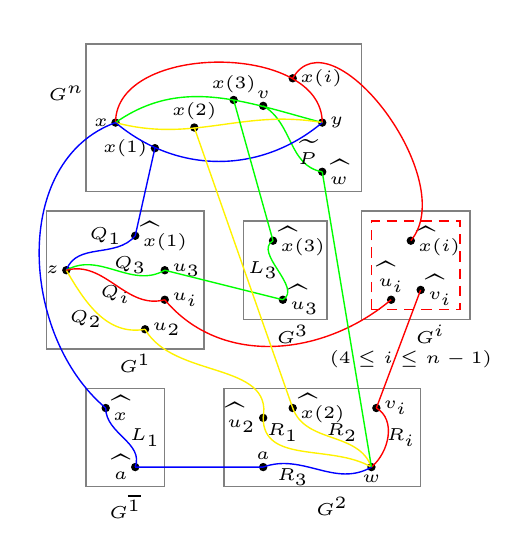}}
\caption{\small Illustration for Subcase 3.3.3} \label{figS211C333}
\end{minipage}

\end{figure}

Let $T_2=P_2\cup R_1\cup Q_2\cup \widetilde{P}(\widehat{w},v) \cup\{u_2\widehat{u}_2, w\widehat{w}\}$, $T_1=P_1\cup L_1\cup Q_1\cup R_2\cup \{x\widehat{x}, x(1)\widehat{x}(1), a\widehat{a}\}$ and $T_i$ ($3\le i\le n-1$) be the same as in Eq.(\ref{EqTi8}). See Figure \ref{figS211C332}. Then $T_1,\cdots, T_{n-1}$ are ($n-1$) IDSTs in $BP_n$.

\vskip 2mm
{\bf Subcase 3.3.3}. $v\in V(P_i)$ for $3\le i\le n-1$.

Without loss of generality, we may assume that $v\in V(P_3)$.

There are ($n-3$) vertices $a, v_4, \cdots, v_{n-1}$ in $V(G^2)\backslash \{\widehat{u}_2, \widehat{x}(2)\}$ such that $\widehat{a}\in V(G^{\overline{1}})$ and $\widehat{v}_i\in V(G^i)$ for $4\le i\le n-1$. Moreover, there is an ($n-1$)-fan $R_1, \cdots, R_{n-1}$ in $G^2$ from $w$ to $\{\widehat{u}_2, \widehat{x}(2), a, v_4, \cdots, v_{n-1}\}$ that $\widehat{u}_2\in V(R_1)$, $\widehat{x}(2)\in V(R_2)$, $a\in V(R_3)$ and $v_i\in V(R_i)$ for $4\le i \le n-1$.

Since $G^i$ is connected for each $i\in [[n]]$, there is a ($\widehat{x}, \widehat{a}$)-path $L_1$ in $G^{\overline{1}}$ and a ($\widehat{x}(3), \widehat{u}_3$)-path $L_3$ in $G^{3}$, respectively.

Let
$$T_3'=P_3\cup \widetilde{P}(\widehat{w},v)\cup L_3\cup Q_3\cup \{x(3)\widehat{x}(3), u_3\widehat{u}_3, w\widehat{w}\},$$
$$T_2=P_2\cup  R_2\cup R_1\cup Q_2\cup \{u_2\widehat{u}_{2}, x(2)\widehat{x}(2)\}$$
and
$$T_1=P_1\cup L_1 \cup Q_1\cup R_3\cup \{x\widehat{x}, x(1)\widehat{x}(1), a\widehat{a}\}.$$
For $4\le i\le n-1$, let $\widehat{T}_i$ and $T_i$ be the same as in Eq.(\ref{EqTi8}). See Figure \ref{figS211C333}. Then $T_1, T_2, T_3', T_4, \cdots, T_{n-1}$ are ($n-1$) IDSTs in $BP_n$.  \hfill$\Box$

\section{$|S\cap V(G^i)|\le 1$ for each $i\in [[n]]$}\label{secS1}

\begin{lemma}\label{lemK4S1}
For $n\ge 3$, let $BP_n=G^1\oplus G^{\overline{1}}\oplus G^2 \oplus G^{\overline{2}}\oplus \cdots \oplus  G^n \oplus G^{\overline{n}}$ and $S=\{x,y,z,w\}$ be any 4-subset of $V(BP_n)$. If $|S\cap V(G^i)|\le 1$ for any integer $i\in[[n]]$. Then there exist ($n-1$) IDSTs in $BP_n$.
\end{lemma}

\noindent{\bf Proof}\; Without loss of generality, we may assume that $x\in V(G^r)$, $y\in V(G^k)$, $z\in V(G^l)$ and $w\in V(G^m)$, where $\{r,k,l,m\}\subseteq [[n]]$.

\vskip 2mm
{\bf Case 1}. $n\ge 7$.

Since $2n-8\ge n-1$ for $n\ge 7$, there exist ($n-1$) different integers $h_1, \cdots, h_{n-1}$ in $[[n]]\backslash\{r,k,l,m,\bar{r},\bar{k},\bar{l},\bar{m}\}$. By Lemma \ref{lemBPn1}, we may choose ($n-1$) vertices $a_i\in V(G^r)$ that $\widehat{a}_i\in V(G^{h_i})$ for $i\in [n-1]$. Analogously, for $i\in [n-1]$, there are ($n-1$) vertices $b_i\in V(G^k)$ that $\widehat{b}_i\in V(G^{h_i})$, there are ($n-1$) vertices $c_i\in V(G^l)$ that $\widehat{c}_i\in V(G^{h_i})$ and ($n-1$) vertices $d_i\in V(G^m)$ that $\widehat{d}_i\in V(G^{h_i})$.

Since $G^{h_i}$ is connected, there is a $\{\widehat{a}_i,\widehat{b}_i,\widehat{c}_i, \widehat{d}_i\}$-tree $\widehat{T}_i$ in $G^{h_i}$ for $i\in [n-1]$.

Let $A=\{a_1,\cdots, a_{n-1}\}$, $B=\{b_1,\cdots, b_{n-1}\}$, $C=\{c_1,\cdots, c_{n-1}\}$ and $D=\{d_1,\cdots, d_{n-1}\}$. Obviously, $|A|=|B|=|C|=|D|=n-1$. By Lemma \ref{lemKfan} and Lemma \ref{lemBPn1}(2), there is an ($n-1$)-fan $P_1, \cdots, P_{n-1}$ in $G^r$ from $x$ to $A$ such that $a_i\in V(P_i)$ for $i\in [n-1]$. Similarly, for $i\in [n-1]$, there is an ($n-1$)-fan $Q_1, \cdots, Q_{n-1}$ in $G^k$ from $y$ to $B$ such that $b_i\in V(Q_i)$, there is an ($n-1$)-fan $R_1, \cdots, R_{n-1}$ in $G^{l}$ from $z$ to $C$ such that
$c_i\in V(R_i)$ and an ($n-1$)-fan $L_1, \cdots, L_{n-1}$ in $G^{m}$ from $w$ to $D$ such that $d_i\in V(L_i)$.

For $i\in [n-1]$, let
$$T_i=\widehat{T}_i\cup P_i\cup Q_i\cup R_i\cup L_i\cup \{a_i\widehat{a}_i, b_i\widehat{b}_i, c_i\widehat{c}_i,d_i\widehat{d}_i\}.$$
Then $T_1, \cdots, T_{n-1}$ are ($n-1$) IDSTs in $BP_n$.

\vskip 2mm
{\bf Case 2}. $n=6$.

\vskip 2mm
{\bf Subcase 2.1}. There exist two integers $i$ and $j$ in $\{r,k,l,m\}$ that $i+j=0$.

Without loss of generality, assume that $r=1$ and $k=\overline{1}$. Since $2n-6>n-1$ for $n=6$, there exist five integers $h_1,\cdots, h_5$ in $[[n]]\backslash\{1,\overline{1}, l, \overline{l}, m, \overline{m}\}$. By similar methods in Case 1, we may construct ($n-1$) IDSTs in $BP_n$.

\vskip 2mm
{\bf Subcase 2.2}. For any two integers $i$ and $j$ in $\{r,k,l,m\}$, it has $i+j\ne 0$.

Without loss of generality, we may assume that $x\in V(G^1)$, $y\in V(G^2)$, $z\in V(G^3)$ and $w\in V(G^4)$.

By Lemma \ref{lemBPn1}, there are five vertices $a_1, \cdots, a_5$ in $V(G^1)$ that $\widehat{a}_1\in V(G^{\overline{2}})$, $\widehat{a}_2\in V(G^{5})$, $\widehat{a}_3\in V(G^{\overline{5}})$, $\widehat{a}_4\in V(G^{6})$ and $\widehat{a}_5\in V(G^{\overline{6}})$. Furthermore, there are different vertices $\{b_1, \cdots, b_5\}\subseteq V(G^2)$, $\{c_1, \cdots, c_5\}\subseteq V(G^3)$ and $\{d_1, \cdots, d_5\}\subseteq V(G^4)$ satisfying  $\{\widehat{b}_1, \widehat{c}_1, \widehat{d}_1\}\subseteq V(G^{\overline{1}})$, $\{\widehat{b}_2, \widehat{c}_2, \widehat{d}_2\}\subseteq V(G^{5})$, $\{\widehat{b}_3, \widehat{c}_3, \widehat{d}_3\}\subseteq V(G^{\overline{5}})$, $\{\widehat{b}_4, \widehat{c}_4, \widehat{d}_4\}\subseteq V(G^{6})$ and $\{\widehat{b}_5, \widehat{c}_5, \widehat{d}_5\}\subseteq V(G^{\overline{6}})$.

By Lemma \ref{lemBPn1}(5), there is a $\{\widehat{a}_1, \widehat{b}_1, \widehat{c}_1, \widehat{d}_1\}$-tree $\widehat{T}_1$ in $G^{\overline{1}}\cup G^{\overline{2}}$. Furthermore, there is a $\{\widehat{a}_2, \widehat{b}_2, \widehat{c}_2, \widehat{d}_2\}$-tree $\widehat{T}_2$ in $G^5$, a $\{\widehat{a}_3, \widehat{b}_3, \widehat{c}_3, \widehat{d}_3\}$-tree $\widehat{T}_3$ in $G^{\overline{5}}$, a $\{\widehat{a}_4, \widehat{b}_4, \widehat{c}_4, \widehat{d}_4\}$-tree $\widehat{T}_4$ in $G^6$ and a $\{\widehat{a}_5, \widehat{b}_5, \widehat{c}_5, \widehat{d}_5\}$-tree $\widehat{T}_5$ in $G^{\overline{6}}$, respectively.

For $1\le i\le 5$, we may define $P_i$, $Q_i$, $R_i$, $L_i$ and $T_i$ be the same as in Case 1. Then $T_1, \cdots, T_{5}$ are five IDSTs in $BP_n$.

\vskip 2mm
{\bf Case 3}. $n=5$.

\vskip 2mm
{\bf Subcase 3.1}. There exist two integers $i$ and $j$ in $\{r,k,l,m\}$ that $i+j=0$.

By similar methods in Subcase 2.1, four IDSTs in $BP_n$ can be constructed and the proof is omitted here.

\vskip 2mm
{\bf Subcase 3.2}. For any two integers $i$ and $j$ in $\{r,k,l,m\}$, it has $i+j\ne 0$.

Without loss of generality, we may assume that $x\in V(G^1)$, $y\in V(G^2)$, $z\in V(G^3)$ and $w\in V(G^4)$.

By Lemma \ref{lemBPn1}, there are four vertices $a_1, a_2, a_3, a_4$ in $V(G^1)$ that $\widehat{a}_1\in V(G^{\overline{2}})$, $\widehat{a}_2\in V(G^{\overline{3}})$, $\widehat{a}_3\in V(G^{5})$ and $\widehat{a}_4\in V(G^{\overline{5}})$. Moreover, there are vertices $\{b_1, \cdots, b_4\}\subseteq V(G^2)$, $\{c_1, \cdots, c_4\}\subseteq V(G^3)$ and $\{d_1, \cdots, d_4\}\subseteq V(G^4)$ satisfying  $\{\widehat{b}_1, \widehat{c}_1, \widehat{d}_1\}\subseteq V(G^{\overline{1}})$, $\{\widehat{b}_2, \widehat{c}_2, \widehat{d}_2\}\subseteq V(G^{\overline{3}})\cup V(G^{\overline{4}})$, $\{\widehat{b}_3, \widehat{c}_3, \widehat{d}_3\}\subseteq V(G^{5})$ and $\{\widehat{b}_4, \widehat{c}_4, \widehat{d}_4\}\subseteq V(G^{\overline{5}})$.

Since $G^{\overline{2i-1}}\cup G^{\overline{2i}}$ is connected, there is a $\{\widehat{a}_i, \widehat{b}_i, \widehat{c}_i, \widehat{d}_i\}$-tree $\widehat{T}_i$ in  $G^{\overline{2i-1}}\cup G^{\overline{2i}}$ for $1\le i\le 2$. Moreover, there is a $\{\widehat{a}_3, \widehat{b}_3, \widehat{c}_3, \widehat{d}_3\}$-tree $\widehat{T}_3$ in $G^{5}$ and a $\{\widehat{a}_4, \widehat{b}_4, \widehat{c}_4, \widehat{d}_4\}$-tree $\widehat{T}_4$ in $G^{\overline{5}}$, respectively.

For $1\le i\le 4$, we may define $P_i$, $Q_i$, $R_i$, $L_i$ and $T_i$ be the same as in Case 1. Then $T_1, \cdots, T_{4}$ are four IDSTs in $BP_n$.

\vskip 2mm
{\bf Case 4}. $n=4$.

\vskip 2mm
{\bf Subcase 4.1}. There exist two pairs of integers in $\{r,k,l,m\}$ that their sum is 0.

By similar arguments in Subcase 2.1, three IDSTs can be constructed and the details are omitted.

\vskip 2mm
{\bf Subcase 4.2}. There are exactly two integers $i$ and $j$ in $\{r,k,l,m\}$ that $i+j=0$.

Without loss of generality, we may assume that $x\in V(G^1)$, $y\in V(G^{\overline{1}})$, $z\in V(G^2)$ and $w\in V(G^3)$.

By Lemma \ref{lemBPn1}, there are different vertices $\{a_1, a_2, a_3\}\subseteq V(G^1)$, $\{b_1, b_2, b_3\}\subseteq V(G^{\overline{1}})$, $\{c_1, c_2, c_3\}\subseteq V(G^2)$ and $\{d_1, d_2, d_3\}\subseteq V(G^3)$ that $\{\widehat{a}_1, \widehat{b}_1, \widehat{c}_1, \widehat{d}_1\}\subseteq V(G^{\overline{2}}) \cup V(G^{\overline{3}})$, $\{\widehat{a}_2, \widehat{b}_2, \widehat{c}_2, \widehat{d}_2\}\subseteq V(G^{4})$ and $\{\widehat{a}_3, \widehat{b}_3, \widehat{c}_3, \widehat{d}_3\}\subseteq V(G^{\overline{4}})$.

According to Lemma \ref{lemBPn1}(5), there is $\{\widehat{a}_1, \widehat{b}_1, \widehat{c}_1, \widehat{d}_1\}$-tree $\widehat{T}_1$ in $G^{\overline{2}}\cup G^{\overline{3}}$ since $G^{\overline{2}}\cup G^{\overline{3}}$ is connected. Moreover, there is a  $\{\widehat{a}_2, \widehat{b}_2, \widehat{c}_2, \widehat{d}_2\}$-tree $\widehat{T}_2$ in $G^4$ and a $\{\widehat{a}_3, \widehat{b}_3, \widehat{c}_3, \widehat{d}_3\}$-tree $\widehat{T}_3$ in $G^{\overline{4}}$, respectively. For $1\le i\le 3$, let $P_i$, $Q_i$, $R_i$, $L_i$ and $T_i$ be the same as in Case 1. Then three IDSTs are constructed.

\vskip 2mm
{\bf Subcase 4.3}. For any two integers $i$ and $j$ in $\{r,k,l,m\}$, it has $i+j\ne 0$.

Without loss of generality, we may assume that $x\in V(G^1)$, $y\in V(G^2)$, $z\in V(G^3)$ and $w\in V(G^4)$.

Based on Lemma \ref{lemBPn2}, there exists an integer $i\in[[4]]\backslash\{1,\overline{1}\}$ that $\widehat{x}\in V(G^i)$ and $\widehat{x}(1)\in V(G^{\overline{i}})$. Without loss of generality, assume that $i=2$. That is to say, $\widehat{x}\in V(G^2)$ and $\widehat{x}(1)\in V(G^{\overline{2}})$. For convenience, denote by $\widehat{x}(1)=g$. Note that $\widehat{g}(1)\in V(G^{\overline{1}})$. There exist two vertices $a_2$ and $a_3$ in $V(G^1)\backslash\{x, x(1)\}$ that $\widehat{a}_i\in V(G^{\overline{i}})$ for $i\in \{2,3\}$. This is possible since $(n-2)!\times 2^{n-2}> 2$ when $n=4$.

There exist two vertices $b_1$ and $b_3$ in $V(G^2)\backslash\{\widehat{x}\}$ that $\widehat{b}_i\in V(G^{\overline{i}})$ for $i\in\{1,3\}$. Moreover, there are three vertices $\{c_1, c_2, c_3\}\subseteq V(G^3)$ that $\widehat{c}_i\in V(G^{\overline{i}})$ for $i\in\{1,2\}$ and $\widehat{c}_3\in V(G^{\overline{4}})$ and three vertices $\{d_1, d_2, d_3\}\subseteq V(G^4)$ that $\widehat{d}_i\in V(G^{\overline{i}})$ where $i\in\{1,2,3\}$.

\begin{figure}[htbp]

\centering
\resizebox{0.98\textwidth}{!} {\includegraphics{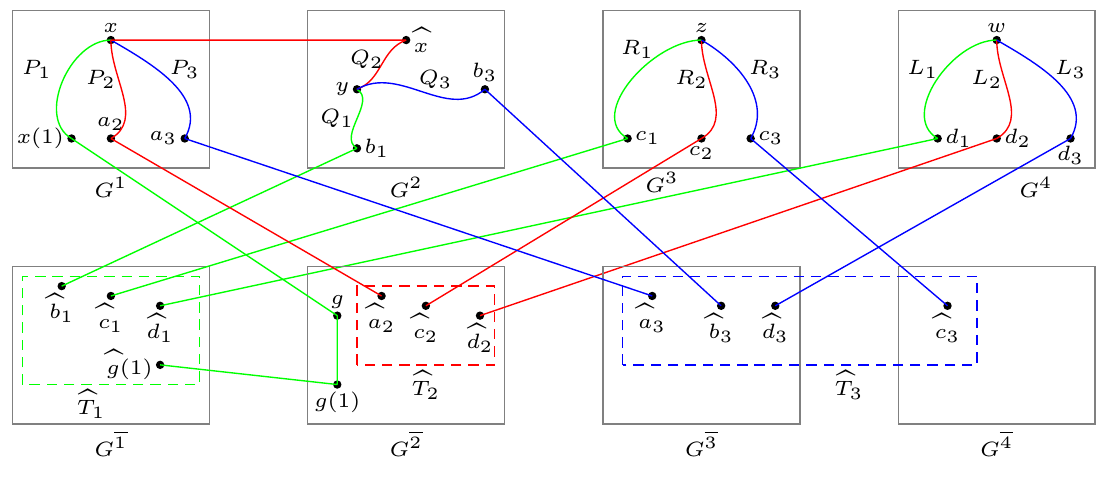}}
\caption{\small Illustration for Subcase 4.3} \label{figS1N4}

\end{figure}

Based on Lemma \ref{lemBPn3}, $G^{\overline{2}}\backslash \{g, g(1)\}$ is still connected. Thus, we can find a $\{\widehat{a}_2, \widehat{c}_2, \widehat{d}_2\}$-tree $\widehat{T}_2$ in $G^{\overline{2}}\backslash \{g, g(1)\}$. Moreover, there is a $\{\widehat{b}_1, \widehat{c}_1, \widehat{d}_1, \widehat{g}(1)\}$-tree $\widehat{T}_1$ in $G^{\overline{1}}$ and a $\{\widehat{a}_3, \widehat{b}_3, \widehat{c}_3, \widehat{d}_3\}$-tree $\widehat{T}_3$ in $G^{\overline{3}}\cup G^{\overline{4}}$.

Remind that $\kappa(G^i)=3$ for $1\le i\le 4$, there is a 3-fan $P_1, P_2$ and $P_3$ in $G^1$ from $x$ to $\{x(1), a_2, a_3\}$ that $x(1)\in V(P_1)$ and $a_i\in V(P_i)$ for $i\in \{2,3\}$. There is a 3-fan $Q_1, Q_2$ and $Q_3$ in $G^2$ from $y$ to $\{\widehat{x}, b_1, b_3\}$ that $\widehat{x}\in V(Q_2)$ and $b_i\in V(Q_i)$ for $i\in \{1,3\}$. There is a 3-fan $R_1, R_2$ and $R_3$ in $G^3$ from $z$ to $\{c_1, c_2, c_3\}$ that $c_i\in V(R_i)$  and a 3-fan $L_1, L_2$ and $L_3$ in $G^4$ from $w$ to $\{d_1, d_2, d_3\}$ that $d_i\in V(L_i)$ for $1\le i\le 3$.

Let
$$T_1=\widehat{T}_1\cup P_1\cup Q_1\cup R_1\cup L_1\cup \{x(1)g, gg(1), g(1)\widehat{g}(1), b_1\widehat{b}_1, c_1\widehat{c}_1, d_1\widehat{d}_1\},$$
$$T_2=\widehat{T}_2\cup P_2\cup Q_2\cup R_2\cup L_2\cup \{x\widehat{x}, a_2\widehat{a}_2, c_2\widehat{c}_2, d_2\widehat{d}_2\}$$
and $T_3=\widehat{T}_3\cup P_3\cup Q_3\cup R_3\cup L_3\cup \{a_3\widehat{a}_3, b_3\widehat{b}_3, c_3\widehat{c}_3, d_3\widehat{d}_3\}.$ See Figure \ref{figS1N4}. It is seen that $T_1, T_2$ and $T_3$ are three IDSTs in $BP_n$.

\vskip 2mm
{\bf Case 5}. $n=3$.

\vskip 2mm
{\bf Subcase 5.1}. There exist two pairs of integers in $\{r,k,l,m\}$ that their sum is 0.

By similar arguments in Subcase 2.1, two IDSTs can be constructed.

\vskip 2mm
{\bf Subcase 5.2}. There are exactly two integers $i$ and $j$ in $\{r,k,l,m\}$ that $i+j=0$.

Without loss of generality, we may assume that $x\in V(G^1)$, $y\in V(G^{\overline{1}})$, $z\in V(G^2)$ and $w\in V(G^3)$.

\begin{figure}[htbp]
\begin{minipage}[t]{0.47\linewidth}
\centering
\resizebox{0.8\textwidth}{!} {\includegraphics{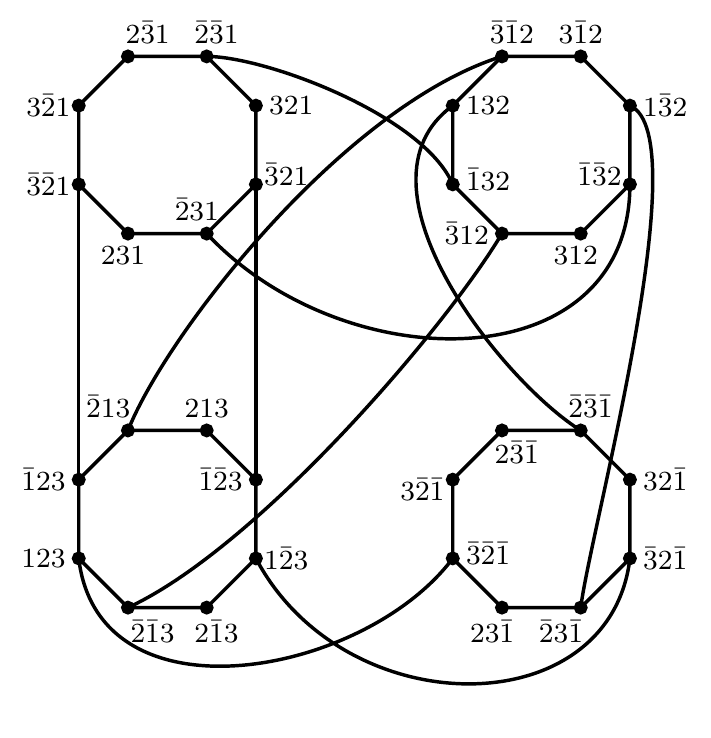}}
\caption{\small The subgraph $H$} \label{figBP3H}
\end{minipage}
\begin{minipage}[t]{0.45\linewidth}
\centering
\resizebox{0.41\textwidth}{!} {\includegraphics{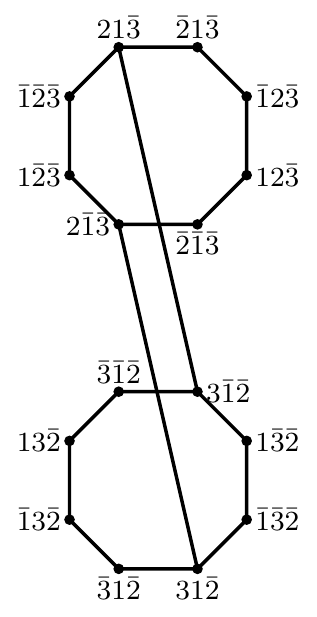}}
\caption{\small The subgraph $H'$} \label{figH1}
\end{minipage}
\end{figure}

Let $H=G^1\cup G^{\overline{1}}\cup G^2\cup G^3$ and $H'=BP_3\backslash H$. See Figure \ref{figBP3H} and Figure \ref{figH1}. The following Claims 1 and 2 can be derived easily.

\vskip 2mm
{\bf Claim 1}. Both $H$ and $H'$ are connected.

\vskip 2mm
{\bf Claim 2}. For $i\in \{\overline{1}, 1, 2, 3\}$, let $a_i$ be any vertex in $V(G^i)$. Then $H\backslash\{a_{\overline{1}}, a_1, a_2, a_3\}$ is connected.

\vskip 2mm
Let $p_1=\bar{3}\bar{1}2$, $p_2=3\bar{1}2$, $p_3=\bar{3}12$, $p_4=312$, $p_5=\bar{2}13$, $p_6=213$, $p_7=\bar{2}\bar{1}3$ and $p_8=2\bar{1}3$ be eight vertices in $V(H)$. Moreover, let $F_1=\{p_1,p_2\}$, $F_2=\{p_3,p_4\}$, $F_3=\{p_5,p_6\}$ and $F_4=\{p_7,p_8\}$ be 2-subsets of $V(H)$.

\vskip 2mm
{\bf Claim 3}. For $i\in \{1, \overline{1}\}$, let $a_i$ be any vertex in $V(G^i)$. Let $V_1$ be a 2-subset of $V(G^2)$ that $V_1\in \{F_1, F_2\}$ and $V_2$ be a 2-subset of $V(G^3)$ that $V_2\in \{F_3,F_4\}$. Then $H\backslash\Big(\{a_{\overline{1}}, a_1\}\cup V_1\cup V_2\Big)$ is connected.

\vskip 2mm
\noindent{\bf Proof of Claim 3}. Note that $d_H(p_{2i})=2$ for $i\in [4]$. By Claim 2, $H\backslash\{a_{\overline{1}}, a_1, a_2, a_3\}$ is connected for $a_2\in\{p_1, p_3\}$, $a_3\in\{p_5, p_7\}$ and $a_i\in V(G^i)$ for $i\in \{1, \overline{1}\}$. It is seen that $H\backslash\Big(\{a_{\overline{1}}, a_1\}\cup V_1\cup V_2\Big)$ is obtained from $H\backslash\{a_{\overline{1}}, a_1, a_2, a_3\}$ by deleting two vertices of degree one. Thus, $H\backslash\Big(\{a_{\overline{1}}, a_1\}\cup V_1\cup V_2\Big)$ is connected. \hfill$\Box$

\vskip 2mm
For $x\in V(G^1)$ and $1\le i\le 2$, we can see that either $\widehat{x}$ or $\widehat{x}(i)$ belong to $V(H')$. Let
\[
P_1(x)= \left\{
\begin{array}{ll}
 x\widehat{x}, \quad & {\rm if}~ \widehat{x}\in V(H'), \\
\{xx(i), x(i)\widehat{x}(i)\}, \quad &  {\rm if}~ \widehat{x}(i)\in V(H').
\end{array}
\right.
\]
We say that $P_1(x)$ is a {\it type-1-path of} $x$. Likewise, there is a type-1-path $P_1(y)$ of $y\in V(G^{\overline{1}})$.

For $z\in V(G^2)$, it is not difficult to check that there is a type-1-path of $z$ if $z\notin \{132, \bar{1}32\}$. If $z=132$, denote by $P_2(z)$ the path $\{zp_1, p_1p_2, p_2\widehat{p}_2\}$. Furthermore, denote by $P_2(z)$ the path $\{zp_3, p_3p_4, p_4\widehat{p}_4\}$ if $z=\bar{1}32$. We say that $P_2(z)$ is a {\it type-2-path of} $z$.

For a vertex $a\in V(H)$ and $1\le i\le 2$, denote by $IP_i(a)$ the set of interior vertices of $P_i(a)$, denote by $TP_i(a)$ the terminal vertex of $P_i(a)$ in $H'$. Note that $|IP_1(a)|\le 1$ for each $a\in V(H)$. In addition, $IP_2(z)\in \{F_1,F_2\}$ if $z\in \{132, \bar{1}32\}$.

Likewise, there is either a type-1-path $P_1(w)$ or a type-2-path $P_2(w)$ of $w\in V(G^3)$ such that $IP_2(w)\in \{F_3,F_4\}.$ Combined with Claims 1, 2 and 3, there is an $S$-tree $T_1$ in $H\backslash \{IP_1(x), IP_1(y), IP_i(z), IP_j(w)\}$ and a $\{TP_1(x), TP_1(y), TP_i(z), TP_j(w)\}$-tree $\widehat{T}_2$ in $H'$, $1\le i\le 2$ and $1\le j\le 2$. Let
$$T_2=\widehat{T}_2\cup P_1(x)\cup P_1(y)\cup P_i(z)\cup P_j(w), \;\; 1\le i\le 2, 1\le j\le 2.$$
Then $T_1$ and $T_2$ are two desired IDSTs in $BP_3$. The proof is done. \hfill$\Box$

\section{The generalized 4-connectivity of the burnt pancake graph}\label{secthm}

Now, we are prepared to prove Theorem \ref{thmBPn}, the main result of our paper.

\vskip 2mm
\noindent{\bf Proof of Theorem \ref{thmBPn}}\; Together with Lemma \ref{lemupperKk} and Lemma \ref{lemBPn1}(1), $\kappa_4(BP_n)\le \delta(BP_n)-1=n-1$ for $n\ge 2$. We shall prove the reverse inequality by induction on $n$. Firstly, $\kappa_4(BP_2)\ge 1$ since $BP_2$ is connected. Now suppose that $n\ge 3$ and the result holds for any integer $m<n$, i.e., $\kappa_4(BP_m)\ge m-1$. Let $S=\{x,y,z,w\}$ be any 4-subset of $V(BP_n)$. The following cases are distinguished.

\vskip 2mm
{\bf Case 1.}\; There exists an integer $i\in [[n]]$ such that $S\subseteq V(G^i)$.

W.l.o.g., we may assume that $\{x, y, z, w\}\subseteq V(G^1)$. By induction hypothesis, there exist ($n-2$) IDSTs $T_1, \cdots, T_{n-2}$ in $G^1$ since $G^1$ is isomorphic to $BP_{n-1}$. Recall that $\widehat{x}$, $\widehat{y}$, $\widehat{z}$ and $\widehat{w}$ are out-neighbours of $x, y, z$ and $w$, respectively. There is a $\{\widehat{x}, \widehat{y},\widehat{z}, \widehat{w}\}$-tree $\widehat{T}_{n-1}$ in $BP_n\setminus V(G^1)$ since $BP_n\setminus V(G^1)$ is connected. Let
$$T_{n-1}=\widehat{T}_{n-1}\cup \{x\widehat{x}, y\widehat{y}, z\widehat{z}, w\widehat{w}\}.$$
Then $T_1, \cdots, T_{n-2}, T_{n-1}$ are ($n-1$) IDSTs in $BP_n$.

\vskip 2mm
{\bf Case 2.}\; There exists an integer $i\in [[n]]$ such that $|S\cap V(G^i)|=3$.

By Lemma \ref{lemK4S3}, ($n-1$) IDSTs can be obtained in $BP_n$.

\vskip 2mm
{\bf Case 3.}\; There exist different integers $\{i,j\} \subseteq [[n]]$ such that $|S\cap V(G^i)|=2$ and $|S\cap V(G^j)|=2$.

According to Lemma \ref{lemK4S22}, ($n-1$) IDSTs can be obtained in $BP_n$.

\vskip 2mm
{\bf Case 4.}\; There are different integers $i,j$ and $l$ in $[[n]]$ such that $|S\cap V(G^i)|=2$ and $|S\cap V(G^j)|=|S\cap V(G^l)|=1$.

Based on Lemma \ref{lemK4S211}, there are ($n-1$) IDSTs in $BP_n$.

\vskip 2mm
{\bf Case 5.}\; For any integer $i\in[[n]]$, $|S\cap V(G^i)|\le 1$.

By Lemma \ref{lemK4S1}, ($n-1$) IDSTs can be obtained in $BP_n$.

In all, there are ($n-1$) IDSTs for any 4-subset $S\subseteq V(BP_n)$. Therefore, $\kappa_4(BP_n)\ge n-1$ for $n\ge 2$. The proof is completed. \hfill$\Box$

\section{Conclusion}\label{seccon}

The generalized $k$-connectivity is a natural generalization of the classical connectivity and can serve for measuring the capability of a network $G$ to connect any $k$ vertices in $G$. In this paper, we focused on the generalized 4-connectivity of the burnt pancake graph $BP_n$. By Lemma \ref{lemupperKk}, the generalized 4-connectivity of an $n$-regular graph is no more than ($n-1$). Thus, our result illustrated that the burnt pancake graph $BP_n$ has a best possible fault-tolerance when the generalized 4-connectivity is applied as the fault-tolerance index. Moreover, $\kappa_3(BP_n)=n-1$ can be derived  directly by Lemma \ref{lemKkk-1} and Theorem \ref{thmBPn}. This result was also investigated in \cite{BPN3}.

%
%
%

%



\end{document}